\documentclass[a4paper,12pt]{article}
\usepackage{amssymb}
\usepackage{amsmath}
\usepackage{amsthm}

\newtheorem{theorem}{Theorem}[section]
\newtheorem{proposition}[theorem]{Proposition}
\newtheorem{lemma}[theorem]{Lemma}
\newtheorem{cor}[theorem]{Corollary}

\newcommand{\A }{\mathbb{A}}
\newcommand{\ad }{\textrm{ad}}

\newcommand{\B }{\mathcal{B}}

\newcommand{\cC}{\mathcal{C}}

\newcommand{\cA }{\mathcal{A}}
\newcommand{\cG}{\mathcal{G}}
\newcommand{\cI}{\mathcal{I}}

\newcommand{\cN}{\mathcal{N}}

\newcommand{\cS}{\mathcal{S}}

\newcommand{\cqg}{\mathbb{C}_q[G]}
\newcommand{\cqgl}{\mathbb{C}_q[G/L_S]}

\newcommand{\C }{\mathbb{C}}
\newcommand{\dif }{\mathrm{d}}
\newcommand{\e }{\mathrm{e}}
\newcommand{\efrak }{\mathfrak{e}}
\newcommand{\End}{\mathrm{End}}

\newcommand{\Fil }{\mathcal{F}}

\newcommand{\GL }{\mbox{GL} }
\newcommand{\gfrak}{\mathfrak{g}}
\newcommand{\gproots}{\overline{R^+_S}}
\newcommand{\gr}{\mathrm{Gr}}
\newcommand{\hfrak}{\mathfrak{h}}

\newcommand{\hght}{\mathrm{ht}}
\newcommand{\id}{\mbox{Id}}

\newcommand{\im}{\mbox{Im}}

\newcommand{\kopr }{\varDelta }
\newcommand{\kow }{\varDelta }

\newcommand{\lfrak}{\mathfrak{l}}

\newcommand{\Lin }{\mathrm{Lin}}
\newcommand{\lord }{\rhd }
\newcommand{\lordne }{\,\raisebox{.5ex}{\makebox[0pt][l]{$\scriptstyle \rhd $}%
\raisebox{-1ex}{$\scriptstyle \nsim $}}\,}

\newcommand{\N }{\mathbb{N}}
\newcommand{\nfrak}{\mathfrak{n}}

\newcommand{\ot }{\otimes }
\newcommand{\pair }[2]{\langle #1,#2 \rangle }
\newcommand{\pfrak}{\mathfrak{p}}
\newcommand{\pord }{\succ }
\newcommand{\pordne }{{\,\scriptstyle \succnsim \,}}

\newcommand{\Proj}{\mathbb{P}}

\newcommand{\Q }{\mathbb{Q}}
\newcommand{\R }{\mathbb{R}}
\newcommand{\rank }{\mbox{rank}}

\newcommand{\rh }{\hat{R}}
\newcommand{\rid }{\mathcal{R}}
\newcommand{\ra }{\acute{R}}
\newcommand{\rg }{\grave{R}}
\newcommand{\rc }{\check{R}}
\newcommand{\rhm }{\hat{R}^-}
\newcommand{\ram }{\acute{R}^-}
\newcommand{\rgm }{\grave{R}^-}
\newcommand{\rcm }{\check{R}^-}

\newcommand{\rlat }{Q}

\newcommand{\sgp}{S[G/P_S]}
\newcommand{\slfrak}{\mathfrak{sl}}
\newcommand{\sofrak}{\mathfrak{so}}
\newcommand{\spfrak}{\mathfrak{sp}}

\newcommand{\sroots }{\pi }

\newcommand{\U }{U}
\newcommand{\ug }{U(\mathfrak{g})}

\newcommand{\uqg }{U_q(\mathfrak{g})}
\newcommand{\Uqbp }{U_q(\mathfrak{b}_+)}
\newcommand{\Uqbm }{U_q(\mathfrak{b}_-)}
\newcommand{\Uqnp }{U_q(\mathfrak{n}_+)}
\newcommand{\Uqnm }{U_q(\mathfrak{n}_-)}
\newcommand{\Ubar }{\overline{\U }}
\newcommand{\wght }{\mathrm{wt}}
\newcommand{\wlat }{P}
\newcommand{\wurz }{\pi }
\newcommand{\vep }{\varepsilon }
\newcommand{\Vtil }{\tilde{V} }
\newcommand{\Z}{\mathbb{Z}}

\title{The Locally Finite Part of the Dual Coalgebra of Quantized Irreducible
  Flag Manifolds}

\author{Istv\'an Heckenberger and Stefan Kolb}

\date{{\footnotesize\textit{Mathematisches Institut, Universit\"at Leipzig,\\
          Augustusplatz 10,
         04109 Leipzig, Germany}\\Istvan.Heckenberger@math.uni-leipzig.de\qquad
       kolb@itp.uni-leipzig.de
      }\\[\baselineskip] June 18, 2003}

\begin{document}

\maketitle

\begin{abstract}
For quantized irreducible flag manifolds the locally finite part of the dual
coalgebra is shown to coincide with a natural quotient coalgebra
$\Ubar$ of $\uqg$. On the way the coradical filtration of $\Ubar$ is
determined. A graded version of the duality between $\Ubar$ and the quantized
coordinate ring is established. This leads to a natural construction of several
examples of quantized vector spaces.

As an application covariant first order differential calculi on quantized
irreducible flag manifolds are classified.

\end{abstract}  

\section{Introduction}
Let $\gfrak$ denote a complex simple Lie algebra and $G$ the
corresponding connected simply connected algebraic group. The
$q$-deformed universal enveloping algebra $\uqg$ can essentially be
recovered from the $q$-deformed coordinate ring $\cqg$ as the dual Hopf
algebra $\cqg^\circ$, \cite{b-Joseph}. Thus $U:=\uqg$ and $\cqg$
constitute dual realizations of the same mathematical object.
It is the main aim of this paper to establish an analogous duality in the
case of quantized irreducible flag manifolds.

Let $P\subset G$ denote a parabolic subgroup with Levi factor $L$. The
$q$-deformed coordinate ring \cite{a-Dk94}, \cite{a-DS992}
\begin{align*}
  \B:=\C_q[G/L]
  :=\{b\in \cqg\,|\, b_{(1)}\,b_{(2)}(k)=\vep(k)b \mbox{ for all } k\in K\}
\end{align*}  
where $K:=U_q(\lfrak)\subset\U$ and $\lfrak$ is the Lie algebra of $L$,
is a right $U$-module algebra. There exists a natural pairing
\begin{align}\label{vorspannpair}
\B\times\Ubar\rightarrow \C
\end{align}
between the right $U$-module algebra $\B$ and the left $U$-module coalgebra
$\Ubar= U/UK^+$ where $K^+=\{k\in K\,|\, \vep(k)=0\}$.
The dual coalgebra $\B^\circ$ generated by the matrix coefficients of all
finite
dimensional representations of $\B$ is a left $U$-module coalgebra. The main
result of the present paper (Theorem \ref{locfin}) is the following
refinement of the duality (\ref{vorspannpair}) for quantized irreducible
flag manifolds:
\begin{align*}
  \B&=\{f\in\Ubar^\ast\,|\, \dim(f U)<\infty\},\\
  \Ubar&=F(\B^\circ,K):=\{f\in\B^\circ\,|\,\dim(K f)<\infty\}
\end{align*}
where $\Ubar^\ast$ denotes the dual vector space of $\Ubar$.
Note that $F(\B^\circ,K)$ can be considered as an analogue of the locally
finite part $F(\cqg^\circ,U)=\{f\in \cqg^\circ\,|\,\dim(\ad(U)f)<\infty\}$
where $\ad$ denotes the adjoint action. Indeed
$(kf)(b)=(k_{(1)}f S(k_{(2)}))(b)$ for all $b\in \B$, $k\in K$ and
$f\in \cqg^\circ$.
The locally finite part $F(\cqg^\circ,U)$ has in principle been determined in
\cite{a-JoLet1} (cp.~also \cite{b-Joseph}, \cite{a-HeckSchm2}).

It has been explained in \cite{a-HK-QHS} that determining the $K$-module
$F(\B^\circ,K)$ is the main step to classify covariant
first order differential calculi over $\B$ in the sense of Woronowicz
\cite{a-Woro2}. As an application of Theorem \ref{locfin} it is
shown that there exist precisely two nonisomorphic irreducible finite
dimensional covariant first order differential calculi over $\B$.

To determine $F(\B^\circ,K)$ several auxiliary results are proven which
are of interest on their own. First, it is necessary to write $\B$ in terms
of generators and relations. It has been shown in \cite{a-DS992},
\cite{a-stok02p} that $\B$ is generated by certain products of matrix
coefficients. Here, a new proof for all generalized flag manifolds is given
relying mainly on the fact that this statement is equivalent to its classical
analogue.

Moreover, in the case of irreducible generalized flag manifolds the coradical
filtration of the connected coalgebra $\Ubar$
is determined. It is proved that the associated graded coalgebra is the graded
dual
of the graded algebra $\bigoplus_{k=0}^\infty(\B^+)^k/(\B^+)^{k+1}$ where
$\B^+=\{b\in\B\,|\,b(1)=0\}$. This yields a natural construction of several
examples of quantized vector spaces such as quantum $n\times m$-matrices or
quantum orthogonal vector spaces \cite{a-FadResTak1}. 

The assertions of this paper concerning the locally finite part
$F(\B^\circ,K)$ and differential calculi over $\B$ are generalizations
of results obtained in the case of quantum complex Grassmann manifolds in
\cite{phd-kolb}.

The ordering of the paper is as follows. Classical generalized flag
manifolds are recalled in Section \ref{genflag} from a point of view also
suitable in the quantum case. In particular the coordinate algebra
$\C[G/L]$ is given in terms of generators and relations.

The third section is devoted to the relevant notions of quantum groups and
quantum flag manifolds. $R$-matrices corresponding to certain representations
of $\U$ appear in the relations of $\C_q[G/L]$. For explicit calculations
it is useful to know how weight spaces are transformed under the action of
these $R$-matrices. Finally $\C_q[G/L]$ is identified with the subalgebra
$\A_\lambda\subset\cqg$ generated by certain products of matrix
coefficients. 

In Section \ref{findimirrgrad} it is shown that $F(\B^\circ,K)$ is connected,
i.e.~the coradical of $F(B^\circ,K)$ is spanned by the counit $\vep\in\Ubar$.
This result is obtained by explicit calculations making consequent use of the
properties of the $R$-matrices.

In Section \ref{U/UK+} the coradical filtration of $\Ubar$ is calculated in
the irreducible case. 
To determine $F(\B^\circ,K)$ in Section \ref{duality} it is shown that the
pairing $\B/(\B^+)^{k+1}\ot C_k\Ubar\rightarrow \C$ is nondegenerate where
$C_k\Ubar$ denotes the elements of degree $k$ with respect to the coradical
filtration. Combined with the results of Section \ref{findimirrgrad}
this implies that $F(\B^\circ,K)=\Ubar$. The graded duality is established and
the resulting coordinate algebras of quantized vector spaces are
discussed.

In the final Section \ref{diffcalc} the notion of covariant first order
differential calculus is recalled and irreducible covariant first order
differential calculi over $\B$ are classified.

Throughout this paper several filtrations are defined in the following way.
Let $A$ denote an algebra generated by the elements of a set $Z$ and $\cS$ a
totally ordered abelian semigroup. Then any map
$\deg:Z\rightarrow \cS$ defines a filtration $\Fil$ of the algebra $A$ as
follows.
An element $a\in A$ belongs to $\Fil_n$, $n\in \cS$, if and only if it can be
written
as a polynomial in the elements of $Z$ such that every occurring summand
$a_{1\dots k}z_{1}\dots z_{k}$, $a_{1\dots k}\in\C$, $z_i\in Z$, satisfies
$\sum_{j=1}^k\deg(z_j)\le n$. Instead of $a\in\Fil_n$ by slight abuse of
notation we will also write $\deg(a)=n$.

For any Hopf algebra $H$ the symbols $\kow$, $\vep$, and $\kappa$ will
denote the coproduct, counit, and antipode, respectively. Sweedler notation
for coproducts $\kow a=a_{(1)}\ot a_{(2)}$, $a\in H$, will be used.
If the antipode $\kappa$ is invertible we will frequently identify left and
right $H$-module structures on a vector space $V$ by $vh=\kappa^{-1}(h)v$,
$v\in V$, $h\in H$.

The authors are very indebted to L.\,L.~Vaksman for his suggestion to study
quantized irreducible flag manifolds.

\section{Generalized Flag Manifolds} \label{genflag}
First, to fix notations some general notions related to Lie algebras are
recalled. Let $\gfrak$ be a finite dimensional complex simple Lie algebra and
$\hfrak\subset \gfrak$ a fixed Cartan subalgebra.
Let $R\subset \hfrak^\ast$ denote the root system
associated with $(\gfrak,\hfrak)$.
Choose an ordered basis $\wurz=\{\alpha_1,\dots,\alpha_r\}$ of simple roots
for $R$. Let $R^+$ and $R^-$ be the set of positive and negative
roots with respect to $\wurz$, respectively.
Moreover, let $\gfrak=\nfrak_+\oplus \hfrak\oplus \nfrak_-$ be the
corresponding triangular decomposition. 
Identify $\hfrak$ with its dual via the Killing form. The induced
non-degenerate symmetric bilinear form on $\hfrak^*$
is denoted by $(\cdot,\cdot)$. The root lattice $Q=\Z R$
is contained in the weight lattice $P=\{\lambda\in\hfrak^\ast\,|\,
(\lambda,\alpha_i)/d_i\in\Z\,\forall \alpha_i\in\wurz\}$ where 
$d_i:=(\alpha_i,\alpha_i)/2$. In order to avoid roots of the deformation
parameter $q$ in the following sections we rescale $(\cdot,\cdot)$ such that
$(\cdot,\cdot):P\times P\rightarrow \Z$.

For $\mu,\nu\in \wlat$ we will write $\mu \pord \nu$ if $\mu-\nu $ is a sum
of positive roots and $\mu\pordne \nu$ if $\mu\pord \nu$ and
$\mu\neq \nu$. As usual we define $Q^+:=\{\mu\in Q\,|\, \mu\pord 0\}$.
The height $\hght:Q_+\rightarrow \N_0$
is given by $\hght(\sum_{i=1}^r n_i \alpha_i)=\sum_{i=1}^r n_i$.
Let $\lordne $ denote the lexicographic ordering on $\Q\ot_\Z\rlat $
with respect to the ordered set $\sroots $ of simple roots: $\mu \lordne 0$,
$\mu \in \Q\ot_\Z\rlat $, if and only if there exists $j\in \{1,\dots,r\}$
such that $\mu =\sum _{i=j}^r\mu _i\alpha _i$, $\mu_i\in \Q$, and $\mu _j>0$.
We write $\mu \lord 0$ if $\mu \lordne 0$ or $\mu =0$. Note in particular
that $\lord$ induces a total ordering on $P$.

The fundamental weights $\omega_i\in\hfrak^\ast$, $i=1,\dots,r$ are
characterized by $(\omega_i,\alpha_j)/d_j=\delta_{ij}$.
Let $P^+$ denote the set of dominant weights, i.\,e.~the $\N_0$-span of
$\{\omega_i\,|\,i=1,\dots,r\}$.
Recall that $(a_{ij}):=(2(\alpha_i,\alpha_j)/(\alpha_i,\alpha_i))$ is the
Cartan matrix of $\gfrak$ with respect to $\wurz$. 

For $\mu\in P^+$  let $V(\mu)$ denote the uniquely determined finite
dimensional irreducible left $\gfrak$-module with highest weight $\mu$.
More explicitly there exists a nontrivial vector $v_\mu\in V(\mu)$ satisfying
\begin{align}
  E v_\mu=0,\quad H v_\mu=\mu(H)v_\mu
  \qquad \mbox{ for all } H\in\hfrak,\, E\in\nfrak_+.
\end{align}  
For any weight vector $v\in V(\mu)$ let $\wght(v)\in P$ denote the weight of
$v$, that is $Hv=\wght(v)(H)v$ for any $H\in \hfrak$.
In particular $\wght(v_1)-\wght(v_2)\in Q$ for all weight vectors
$v_1,v_2\in V(\mu)$.

Let $W$ denote the Weyl group of $\gfrak$. Recall that for any finite
dimensional $\gfrak$-module $V$ the Weyl group permutes the weight spaces
$V_\mu$ of $V$ and $\dim V_{w\mu}=\dim V_\mu$ for all $w\in W$. If $w_0$ is
a longest element in $W$ then $w_0\mu$ is the lowest weight of the
irreducible representation $V(\mu)$.

Let $G$ denote the connected, simply connected complex Lie group with Lie
algebra $\gfrak$. Recall that $G$ is complex algebraic and let $\C[G]$ denote
its coordinate ring. The finite dimensional representations of $\gfrak$ are
in one-to-one correspondence to the finite dimensional rational $G$-modules.
For $v\in V(\mu)$, $f\in V(\mu)^\ast$ the matrix coefficient
$c^\mu_{f,v}\in \C[G]$ is defined by
\begin{align*}
  c^\mu_{f,v}(g)=f(gv).
\end{align*}
The linear span of matrix coefficients of $V(\mu)$
\begin{align}
  C^{V(\mu)}:=\Lin_\C\{c^\mu_{f,v}\,|\,v\in V(\mu), f\in V(\mu)^\ast\}
\end{align}
obtains a $G$-bimodule structure by
\begin{align}
  (hc^\mu_{f,v}k)(g)=f(kghv)=c^\mu_{fk,hv}(g),\qquad g,h,k\in G.
\end{align}
Here $V(\mu)^\ast$ is considered as a right $G$-module. In the same way
$C^{V(\mu)}$ can be endowed with a $U(\gfrak)$-bimodule structure.
As $G$ is a closed subgroup of $\GL(n)$ for some $n\in \N$ one obtains
\begin{align}\label{peterw}
  \C[G]\cong\bigoplus_{\mu\in P^+} C^{V(\mu)}.
\end{align}

For any set $S\subset \wurz$ of simple roots define $R_S^\pm:=\Z S\cap R^\pm$
and $\overline{R_S^\pm}:=R^\pm\setminus R_S^\pm$. 
Let $P_S$ and $P_S^{op}$ denote the corresponding standard parabolic
subgroups of $G$ with Lie algebra
\begin{align}
  \pfrak_S=\hfrak\oplus\bigoplus_{\alpha\in R^+\cup R^-_S}\gfrak_\alpha,\qquad
  \pfrak^{op}_S=\hfrak\oplus\bigoplus_{\alpha\in R^-\cup R_S^+}\gfrak_\alpha
\end{align}
The generalized flag manifold $G/P_S$ is called irreducible if the
adjoint representation of $\pfrak_S$ on $\gfrak/\pfrak_S$ is irreducible.
Equivalently, $S=\pi\setminus\{\alpha_i\}$ where $\alpha_i$ appears in any
positive root with coefficient at most one. For a complete list of all
irreducible flag manifolds consult e.g.~\cite[p.~27]{b-BastonEastwood}.

One can associate several algebras of functions to the generalized flag
manifold $G/P_S$. We identify the fixed subset
$S=\{\alpha_{i_1},\cdots,\alpha_{i_s}\}$ with
the corresponding subset $\{i_1,\dots,i_s\}\subset \{1,\dots,r\}$.

First consider the irreducible representation $V(\lambda)$ of $\gfrak$ of
highest weight $\lambda=\sum_{i\in\wurz\setminus S}\omega_i$.
Then $G/P_S$ is isomorphic
to the $G$-orbit of the highest weight vector $v_\lambda\in V(\lambda)$ in
projective space $\mathbb{P}(V(\lambda))$. Therefore the homogeneous
coordinate ring $\sgp$ of $G/P_S$ coincides with the subalgebra of $\C[G]$
generated by the matrix coefficients
$\{c^\lambda_{f,v_\lambda}\,|\, f\in V(\lambda)^\ast\}$.

Recall that there is an isomorphism of right $\ug$-module algebras
\begin{align*}
  S[G/P_S]\overset{\cong}{\longrightarrow}\bigoplus_{n=0}^\infty V(n\lambda)^\ast
\end{align*}
where the multiplicative
structure on the right hand side is given by the Cartan multiplication
\begin{align*}
  V(n_1\lambda)^\ast \ot V(n_2\lambda)^\ast
  \rightarrow V((n_1+n_2)\lambda)^\ast.
\end{align*}  
Moreover by a theorem of Kostant (cf.~\cite{b-FultonHarris}) the algebra
$S[G/P_S]$
is quadratic, i.e.~
\begin{align}\label{S[G/P]quadratic}
  S[G/P_S]\cong \bigoplus_{n=0}^\infty (V(\lambda)^\ast)^{\ot n}/\cI(\lambda),
\end{align}  
where $\cI(\lambda)$ denotes the ideal in the tensor algebra generated by the
subspace
\begin{align*}
  \bigoplus_{\makebox[0cm]{$\mu\neq 2\lambda\atop
  V(\mu)\subset V(\lambda)\ot V(\lambda)$}}V(\mu)^\ast\subset
  (V(\lambda)^\ast)^{\ot 2}.
\end{align*}

Similarly  $G/P_S^{op}$ is isomorphic to the $G$-orbit of the lowest weight
vector $f_{-\lambda}\in V(-w_0\lambda)\cong V(\lambda)^*$ in
$\Proj(V(\lambda)^*)$. Therefore $S[G/P^{op}_S]$ is isomorphic to the
subalgebra of $\C[G]$ generated by the matrix coefficients
$\{c^{-w_0\lambda}_{v,f_{-\lambda}}\,|\,
v\in V(-w_0\lambda)^\ast\cong V(\lambda)\}$. Again there exists an
isomorphism
\begin{align*}
S[G/P_S^{op}]\overset{\cong}{\longrightarrow}
\bigoplus_{n=0}^\infty V(n\lambda)
\end{align*}
of right $U(\gfrak)$-module algebras and $S[G/P_S^{op}]$ is
quadratic.

Let $\A_\lambda \subset \C[G]$ denote the subalgebra generated by the elements
\begin{align}\label{Alambdagen}
  \left\{z_{fv}:=\frac{c^\lambda_{f,v_\lambda} c^{-w_0\lambda}_{v,f_{-\lambda}}}
  {f_{-\lambda}(v_\lambda)}
  \,\Bigg|\,
            f\in V(\lambda)^\ast,\, v\in V(\lambda)\right\}.
\end{align}
By construction the space of $n$-fold products of the generators $z_{fv}$
of $\A_\lambda$ is isomorphic to $V(n\lambda)^\ast\ot V(n\lambda)$.
Define $N:=\dim V(\lambda)$ and $I:=\{1,\dots,N\}$.
If $\{v_i\,|\,i\in I\}$ and $\{f_i\,|\, i\in I\}$ are dual bases of
$V(\lambda)$ and $V(\lambda)^\ast$, respectively, then
\begin{align*}
  \sum_{i\in I} z_{f_iv_i}=f_{-\lambda}(v_\lambda)^{-1}\sum_{i\in I}
c^\lambda_{f_i,v_\lambda} c^{-w_0\lambda}_{v_i,f_{-\lambda}}=1.
\end{align*}
Therefore there is a natural inclusion of the space of $n$-fold products of the
generators $z_{fv}$ of $\A_\lambda$ into the space of $(n+1)$-fold products
given by multiplication with $\sum_{i\in I}z_{f_iv_i}$.
With respect to these inclusions $\A_\lambda$ can be written as a direct
limit
\begin{align}\label{Alim}
  \A_\lambda\cong \underrightarrow{\lim} V(n\lambda)^\ast\ot V(n\lambda).
\end{align}  
Thus a complete set of defining relations of $\A_\lambda$ is given by
linearity in the first an second index of $z$ and
\begin{align}
  \sum_i z_{g_iu}z_{h_iw}&=0\qquad \mbox{ if }\sum_i g_i\ot
  h_i\in\bigoplus_{\mu\neq 2\lambda} V(\mu)^\ast,\label{relations1}\\
 \sum_i z_{gu_i}z_{hw_i}&=0\qquad \mbox{ if }\sum_i u_i\ot
  w_i\in\bigoplus_{\mu\neq 2\lambda} V(\mu),\label{relations2}\\
 \sum_{i\in I} z_{f_iv_i}&=1,\label{relations3}
\end{align}  
where as above $\{v_i\,|\,i\in I\}$ is a basis of $V(\lambda)$ with dual basis
$\{f_i\,|\, i\in I\}$.

The algebra $\A_\lambda$ has a straightforward geometric interpretation.
Consider the Levi factor
\begin{align}
\lfrak_S=\hfrak\oplus\bigoplus_{\alpha\in \Z S} \gfrak_\alpha
\end{align}
of $\pfrak_S$ and let $L_S=P_S\cap P_s^{op}\subset G$ denote the
corresponding subgroup.
As $L_S$ is reductive the quotient $G/L_S$ is an affine algebraic
variety with coordinate ring \cite{a-HoKo62}
\begin{align*}
  \C[G/L_S]=\C[G]^{L_S}:=
  \{a\in\C[G]\,|\,a(g)=a(gl)\, \forall g\in G,\,l\in L_S\}.
\end{align*}
Note that by construction $\A_\lambda \subset \C[G/L_S]$.
It is our next aim to show that $\A_\lambda= \C[G/L_S]$.

Being reduced $\A_\lambda$ is the coordinate ring of an affine algebraic
variety. The corresponding algebraic set $Z(\A_\lambda)$ can be identified
with
\begin{align}
  \begin{aligned}
 \{z\in \End(V(\lambda))\,|\, \mbox{trace}(z)=1,\,\rank\, z\le 1,\,
                   \im z\in G v_\lambda\subset \Proj(V(\lambda)),&\\
                   \im z^{\mbox{\tiny t}}\in G f_{-\lambda}\subset
                   \Proj(V(\lambda)^\ast)\}&
  \end{aligned}                 
\end{align}
where $z^{\mbox{\tiny t}}$ denotes the transposed map.
Indeed, $\rank\, z\le 1$ is equivalent to (\ref{relations1}) for antisymmetric
elements $\sum_i g_i\ot h_i$. Moreover, in view of (\ref{S[G/P]quadratic})
Equation (\ref{relations1}) with $u=v$ also implies that
$\im z\in G v_\lambda\subset \Proj(V(\lambda))$.

One now verifies that the map
\begin{align}\label{ZAlam}
  Z(\A_\lambda)&\rightarrow \{(v,f)\in Gv_\lambda\times Gf_{-\lambda}\subset
  \Proj(V(\lambda))\times \Proj(V(\lambda)^\ast)\,|\, f(v)\neq 0\}\\
  z&\mapsto (\im z,\im z^{\mbox{\tiny t}})\nonumber
\end{align}
is an isomorphism of quasiprojective varieties with inverse morphism
\begin{align}\label{inverse}
  (v,f)\mapsto z=(z_{ij})_{i,j=1,\dots,N},\qquad
  z_{ij}=\frac{f_i(v)v_j(f)}{f(v)}
\end{align}  
where $z\in \End(V(\lambda))$ is given with respect to the basis chosen
above. 
Using the identification (\ref{ZAlam}) one obtains a morphism
\begin{align*}
  \psi:G/L_S\rightarrow  Z(\A_\lambda),\qquad
  g\mapsto (gv_\lambda,g f_{-\lambda}).
\end{align*}
Note that by (\ref{inverse}) this morphism corresponds to the inclusion
$\A_\lambda \hookrightarrow \C[G/L_S]$ of coordinate rings. 
\begin{proposition}\label{A=CGL}
  The map $\psi:G/L_S\rightarrow Z(\A_\lambda)$ is an isomorphism of affine
  algebraic varieties. In particular $\A_\lambda\hookrightarrow \C[G/L_S]$
  is an isomorphism.
\end{proposition}  
\begin{proof}
  Let $W_S\subset W$ denote the Weyl group of the Levi factor $L_S$ of $P_S$
  and $W^S:=W/W_S$. Recall the Bruhat decomposition
  (cf. e.\,g.~\cite{b-BastonEastwood})
  \begin{align}\label{bruhat}
    G=\coprod_{w\in W^{-w_0S}} B^{op}wP^{op}_{-w_0S}
  \end{align}
  of $G$ with respect to the parabolic subgroup $P^{op}_{-w_0S}\supset
  B^{op}$ of $G$.
  To verify surjectivity of $\psi$ consider any $(ghv_\lambda,gf_{-\lambda})\in
  G v_\lambda\times Gf_{-\lambda}$, $g,h\in G$, such that $f_{-\lambda}(hv_\lambda)\neq
  0$. Then $hv_\lambda=h w_0v_{w_0\lambda}$ for a lowest weight vector
  $v_{w_0\lambda}\in V(\lambda)$ and hence (\ref{bruhat}) implies that
  \begin{align*}
    h w_0=bw_0p\in B^{op}w_0P^{op}_{-w_0S}.
  \end{align*}
  Therefore $hv_\lambda=bw_0v_{w_0\lambda}\in B^{op}v_\lambda\subset
  \Proj(V(\lambda))$ and hence
  \begin{align*}
    (ghv_\lambda,gf_{-\lambda})=(g bv_\lambda,g f_{-\lambda})=
    (g bv_\lambda,g b f_{-\lambda})\in
     \Proj(V(\lambda))\times \Proj(V(\lambda)^\ast)
  \end{align*}
  which implies surjectivity of $\psi$.

  Since $L_S=P_S\cap P_S^{op}$ the morphism $\psi$ is a bijection and by
  \cite[Thm.~4.6]{b-Humphreys95}
  also birational. Now \cite[Prop.~4.7]{b-Humphreys95} implies that there
  exists a nonempty open set $U\subset  Z(\A_\lambda)$ such that $\psi$
  induces an isomorphism of $\psi^{-1}(U)$ onto $U$. As $\psi$ is compatible
  with the transitive action of $G$ on $G/L_S$ and $Z(\A_\lambda)$ this implies
  that $\psi$ is an isomorphism of affine varieties. 
\end{proof}

\section{Quantum Groups and Quantum Flag Manifolds}\label{qflag}
We keep the notations of the previous section. Let $0\neq q\in \C$ be not a
root of unity. The $q$-deformed universal enveloping algebra $\U=\uqg$
associated
to $\gfrak$
can be defined to be the complex algebra with generators $K_i,K_i^{-1},
E_i,F_i$, $i=1,\dots,r$, and relations 
\begin{align}
\begin{aligned}
&\begin{aligned}
K_iK_i^{-1}&=K_i^{-1}K_i=1,& K_iK_j&=K_jK_i,\\
K_iE_j&=q^{(\alpha_i,\alpha_j)}E_jK_i,& K_iF_j&=q^{-(\alpha_i,\alpha_j)}F_jK_i,
\end{aligned}&\\
&\begin{aligned}
E_iF_j-F_jE_i&=\delta_{ij}\frac{K_i-K_i^{-1}}{q_i-q_i^{-1}},\\
\sum_{k=0}^{1-a_{ij}}(-1)^k\left(\begin{array}{c}1-a_{ij}\\k \end{array}
        \right)_{q_i}& E_i^{1-a_{ij}-k}E_jE_i^k=0,&i&\neq j,\\
\sum_{k=0}^{1-a_{ij}}(-1)^k\left(\begin{array}{c}1-a_{ij}\\k \end{array}
        \right)_{q_i}&F_i^{1-a_{ij}-k}F_jF_i^k=0,&i&\neq j,
\end{aligned}&
\end{aligned}
\end{align}
where $q_i:=q^{d_i}$ and the $q$-deformed binomial coefficients are defined
by
\begin{align*}
{n\choose k}_q=\frac{[n]_q[n{-}1]_q\dots[n{-}k{+}1]_q}{[1]_q[2]_q\dots[k]_q},
    \qquad\mbox{where } [x]_q=\frac{q^x-q^{-x}}{q-q^{-1}}.
\end{align*}
The algebra $\U$ obtains a Hopf algebra structure by
\begin{align}\label{hopfstruc}
\begin{aligned}
\kopr K_i&= K_i\otimes K_i,&\kopr E_i&= E_i\otimes K_i+1\otimes E_i,&
\kopr F_i&= F_i\otimes 1 + K_i^{-1}\otimes F_i,\\
\epsilon(K_i)&=1,&\epsilon(E_i)&=0,& \epsilon(F_i)&=0,\\
\kappa(K_i)&=K_i^{-1},&\kappa(E_i)&=-E_iK_i^{-1},& \kappa(F_i)&=-K_iF_i.
\end{aligned}
\end{align}
Let $\Uqnp,\Uqbp,\Uqnm,\Uqbm\subset \U$ denote the subalgebras generated
by $\{E_i\,|\,i=1,\dots,r\}$, $\{E_i,K_i,K_i^{-1}\,|\,i=1,\dots,r\}$,
$\{F_i\,|\,i=1,\dots,r\}$, and $\{F_i,K_i,K_i^{-1}\,|\,i=1,\dots,r\}$,
respectively.

For $\mu\in P^+$  let $V(\mu)$ denote the uniquely determined finite
dimensional irreducible left $\U$-module with highest weight $\mu$.
More explicitly there exists a nontrivial highest weight vector
$v_\mu\in V(\mu)$ satisfying
\begin{align}
  E_iv_\mu=0,\qquad
  K_iv_\mu=q^{(\mu,\alpha_i)}v_\mu \qquad
  \quad \mbox{ for all } i=1,\dots,r.
\end{align}
A finite dimensional $\U$-module $V$ is called \textit{ of type 1} if
$V\cong \bigoplus_i  V(\mu_i)$ is isomorphic to a direct sum of finitely many
$V(\mu_i)$, $\mu_i\in P^+$. The category $\cC$ of $\U$-modules of type 1
is a tensor category. By this we mean that $\cC$ contains the trivial
$\U$-module $V(0)$ and satisfies
\begin{align}\label{tenscat}
  X,Y\in\mathcal{C} \Rightarrow X\oplus Y,\,X\ot Y,\,X^*\in \mathcal{C}
\end{align}
where $(uf)(x):=f(\kappa(u)x)$ for all $u\in U$, $f\in X^\ast$, $x\in X$.

Moreover, $\cC$ is a braided tensor category \cite{b-CP94}, \cite{b-KS}.
For all $V,W\in \cC$ the braiding
\begin{align*}
  \rh_{V,W}:V\ot W\rightarrow W\ot V
\end{align*}
satisfies
\begin{align}
  \rh_{V,W}(v\ot w)&=q^{(\wght (v),\wght (w))}w\ot v +\label{Rstruk}\\
  &+\sum_{j=1}^r
  q^{(\wght (v)+\alpha_j,\wght (w)-\alpha_j)}(q^{d_j}-q^{-d_j})F_jw\ot E_j v +
  \sum _iw_i\ot v_i\nonumber
\end{align}
where $\wght (v_i)\pordne\wght (v)$, $\wght (w)\pordne\wght (w_i)$, and
$\hght(\wght (v_i)-\wght (v))\ge 2$, $\hght(\wght (w)-\wght (w_i))\ge 2$.  
To simplify notation we will also write
$\rh_{\mu,\nu}:=\rh_{V(\mu),V(\nu)}$ if $\mu,\nu\in P^+$.

To write the $q$-analogues of $\C[G/L_S]$ in terms of generators and
relations similar to (\ref{relations1})--(\ref{relations3}) it will be
helpful to introduce
additional notations for certain special cases of $\rh$.
Recall that $\lambda=\sum_{i\notin S}\omega_i$,
$N=\dim V(\lambda)$, and $I=\{1,\dots,N\}$. Choose a basis
$\{v_i\,|\,i \in I\}$ of weight vectors of $V(\lambda)$ and
let $\{f_i\,|\,i\in I\}$ be the corresponding dual basis.
Define matrices $\rh$, $\rc$, $\ra^-$, and $\rg^-$ by
\begin{align*}
  \rh_{\lambda,\lambda}(v_i{\ot} v_j)&=:
                     \sum_{k,l\in I}\rh^{kl}_{ij} v_k{\ot} v_l,&
  \rh_{-w_0\lambda,-w_0\lambda}(f_i{\ot} f_j)&=:
                     \sum_{k,l\in I}\rc^{kl}_{ij} f_k{\ot} f_l,\\
  \rh_{-w_0\lambda,\lambda}(f_i{\ot} v_j)&=:
                     \sum_{k,l\in I}\ram {}^{kl}_{ij} v_k{\ot} f_l,&
  \rh_{\lambda,-w_0\lambda}(v_i{\ot} f_j)&=:\sum_{k,l\in I}
  \rgm {}^{kl}_{ij} f_k{\ot} v_l.
\end{align*}  
Alternatively
\begin{align*}
  (f_i{\ot} f_j){\circ} \rh_{\lambda,\lambda}&=
                       \sum_{k,l\in I}\rh_{kl}^{ij} f_k{\ot} f_l,&
  (v_i{\ot} v_j){\circ} \rh_{-w_0\lambda,-w_0\lambda}&=\sum_{k,l\in I}
                       \rc_{kl}^{ij} v_k{\ot} v_l,\\
  (f_i{\ot} v_j){\circ} \rh_{-w_0\lambda,\lambda}&=
                       \sum_{k,l\in I}\ram {}_{kl}^{ij} v_k{\ot} f_l,&
  (v_i{\ot} f_j){\circ} \rh_{\lambda,-w_0\lambda}&=
                       \sum_{k,l\in I}\rgm {}_{kl}^{ij} f_k{\ot} v_l,\\
\end{align*}  
where the elements of $V(\lambda)$ are considered as functionals on
$V(\lambda)^\ast$. Let $\rh^-$, $\rc^-$, $\rg$, and $\ra$ denote the
inverse of the matrix $\rh$, $\rc$, $\ra^-$, and $\rg^-$, respectively.

By (\ref{Rstruk}) the matrix $\rh $ has the property that
$\rh ^{ij}_{kl}\not=0$
implies that $i=l,j=k$ or both $\wght (v_j)\pordne\wght (v_k)$ and
$\wght (v_l)\pordne\wght (v_i)$. Therefore we associate to $\rh $ the symbol
$<$ which denotes the positions of the larger weights.
Similar properties are fulfilled for the other types of $R$-matrices.
For example, the relation $\ram {}^{ij}_{kl}\neq 0$ implies that
$i=l,j=k$ or both
$\wght (v_k)\pordne\wght (v_j)$ and $\wght (v_l)\pordne\wght (v_i)$.
We collect these properties in the following table.

\begin{gather}\label{rtabelle}
\begin{array}{|c|c|c|c|c|c|c|c|}
\hline
\rh \rule{0pt}{3ex} & \rhm & \ra & \ram & \rc & \rcm & \rg & \rgm \\
\hline
< & > & \vee & \wedge & > & < & \wedge & \vee \\
\hline
\end{array}
\end{gather}
For $j\in \{1,2,\ldots ,r\}$ let $U_j$ denote the Hopf subalgebra of $\U $
generated by $E_j,F_j,K_l,K_l^{-1}$, $l=1,2,\ldots ,r$.
Moreover, let  $V(\lambda )=\bigoplus _m V(\lambda )_m$ denote the
decomposition of $V(\lambda)$ into irreducible $U_j$-modules.
The following Lemma will be used only in Step 5 of the proof of Proposition
\ref{IrrGradRep}.

\begin{lemma}\label{l-divers}
Let $\{v_i\,|\,i\in I\}$ be a weight basis of $V(\lambda)$ respecting
the decomposition
$\bigoplus _m V(\lambda )_m$ and let $\{f_i\,|\,i\in I\}$ denote the dual
basis.
Assume that there exist $i,k$ such that
$F_jv_i=0$ and $E_jv_i=v_k$. Set $\mu =\wght (v_i)$. Then the following
relations hold:

\begin{itemize}
\item[(i)] \begin{minipage}[t]{12cm}$F_jv_k=
(q^{-(\mu,\alpha _j)}-q^{(\mu,\alpha _j)})/(q^{d_j}-q^{-d_j})v_i$,\quad
$E_jf_k=-q^{-(\mu,\alpha _j)}f_i$,

\noindent $K_jf_i=q^{-(\mu,\alpha _j)}f_i$,\quad $E_jf_i=0$,\quad
$F_jf_i=(q^{2(\mu,\alpha _j)}-1)/(q^{d_j}-q^{-d_j})f_k$.
\end{minipage}

\item[(ii)]
  $\rgm {}^{kk}_{ii}=q^{-(\mu+\alpha_j,\mu+\alpha_j)}(q^{2(\mu,\alpha_j)}-1)$,
  $\rh^{ik}_{ia}=-\delta_{ka}q^{(\mu,\mu)}(q^{2(\mu,\alpha_j)}-1)$.

\item[(iii)]
  $\ra {}^{kk}_{ai}=-\delta_{ai}q^{(\mu ,\mu )}(q^{2(\mu,\alpha _j)}-1)$,
$\rhm {}^{ki}_{ai}=\delta _{ka} q^{-(\mu ,\mu )}(1-q^{-2(\mu ,\alpha _j)})$.

\end{itemize}
\end{lemma}

\begin{proof}
(i) follows from $xf_m(v_n)=f_m(\kappa(x)v_n)$ for all $x\in \U $ and the fact
that the basis $\{v_i\,|\,i\in I\}$ respects the decomposition
$V(\lambda)=\bigoplus _m V(\lambda )_m$ into irreducible $U_j$-modules.

(ii) follows from (i) and Equation (\ref{Rstruk}) for $\rgm$ and $\rh$.

(iii) One has
$0=\delta_{ka}\delta_{ki}=\sum _{l,m\in I}\rgm {}^{kk}_{lm}\ra {}^{lm}_{ai}$.
Since $v_k=E_jv_i$ it follows from table (\ref{rtabelle}) that only
summands with $(l,m)=(i,i)$ or $(l,m)=(k,k)$ can be nonzero.
Thus $\ra {}^{kk}_{ai}=-q^{(\mu +\alpha _j,\mu +\alpha _j)+(\mu ,\mu )}
\delta_{ai}\rgm {}^{kk}_{ii}=-q^{(\mu ,\mu )}(q^{2(\mu,\alpha _j)}-1)
\delta_{ai}$ by (ii).
The second formula follows similarly from $(\rhm {}\rh )^{ki}_{ia}=0$.
\end{proof}

The $q$-deformed coordinate ring $\cqg$ is defined to be the subspace
of the linear dual $\U^\ast$ spanned by the matrix coefficients of the
finite dimensional irreducible representations $V(\mu)$, $\mu\in
P^+$.
For $v\in V(\mu)$, $f\in V(\mu)^\ast$ the matrix coefficient
$c^\mu_{f,v}\in \U^\ast$ is defined by
\begin{align*}
  c^\mu_{f,v}(X)=f(Xv).
\end{align*}
The linear span of the matrix coefficients of $V(\mu)$
\begin{align}
  C^{V(\mu)}:=\Lin_\C\{c^\mu_{f,v}\,|\,v\in V(\mu), f\in V(\mu)^\ast\}
\end{align}
obtains a $\U$-bimodule structure by
\begin{align}
  (Yc^\mu_{f,v}Z)(X)=f(ZXYv)=c^\mu_{fZ,Yv}(X).
\end{align}
Here $V(\mu)^\ast$ is considered as a right $\U$-module. Note that by
construction
\begin{align}\label{qpeterw}
  \cqg\cong\bigoplus_{\mu\in P^+} C^{V(\mu)}
\end{align}
is a Hopf algebra and the pairing
\begin{align}\label{qpair}
  \cqg\ot \U\rightarrow \C
\end{align}
is nondegenerate.

The algebras $S[G/P_S]$, $\A_\lambda$ and $\C[G/L_S]$ have natural
analogues in the $q$-deformed setting. One defines $S_q[G/P_S]\subset \C_q[G]$ 
as the subalgebra generated by the matrix coefficients
$\{c^\lambda_{f,v_\lambda}\,|\, f\in V(\lambda)^\ast\}$
\cite{b-CP94}, \cite{a-LaksResh92}, \cite{a-TaftTo91}, \cite{a-Soib92}.
Again $S_q[G/P_S]\cong\bigoplus_{n=0}^\infty V(n\lambda)^\ast$ endowed with
the Cartan multiplication and
\begin{align*}
  S_q[G/P_S]\cong \bigoplus_{n=0}^\infty
          (V(\lambda)^\ast)^{\ot n}/\cI(\lambda),
\end{align*}  
where the $\U$-module $\cI(\lambda)$ is defined as in Section \ref{genflag}
(\cite{a-TaftTo91}, \cite{a-Brav94}).
The matrix coefficients $\{c^\lambda_{f_i,v_\lambda}\,|\,i\in I\}$ satisfy
the relations
\begin{align}\label{c-lamlam}
  c_{f_i,v_{\lambda}}^{\lambda}c_{f_j,v_\lambda}^\lambda=
  q^{-(\lambda,\lambda)}\sum_{k,l\in I}\rh_{kl}^{ij}c_{f_k,v_\lambda}^\lambda
  c_{f_l,v_{\lambda}}^{\lambda}
\end{align}
As the eigenvalues of $\rh_{\lambda,\lambda}$ are different from
$q^{(\lambda,\lambda)}$ on all subspaces
$V(\mu)\subset V(\lambda)\ot V(\lambda)$, $\mu\neq 2\lambda$, the relations
(\ref{c-lamlam}) form a complete set of defining relations of $S_q[G/P_S]$.

Let $\A_\lambda^q\subset \C_q[G]$ denote the subalgebra generated by the
elements $z_{fv}$ as in (\ref{Alambdagen}). To shorten notation define
$z_{ij}=z_{f_iv_j}$ with $v_j$ and $f_i$ as above. It follows from
(\ref{c-lamlam}) and
\begin{align}\label{c-lamclam}
  c_{v_i,f_{-\lambda}}^{-w_0\lambda}c_{f_j,v_\lambda}^\lambda=
  q^{(\lambda,\lambda)}\sum_{k,l\in I}(\rg^-)_{kl}^{ij}
  c_{f_k,v_\lambda}^\lambda c_{v_l,f_{-\lambda}}^{-w_0\lambda}
\end{align}
that the following relations hold in $\A_\lambda^q$:
\begin{align}
  \sum_{m,n,p,t\in I}\rh^{ij}_{nm}\ra^{mk}_{pt}z_{np}z_{tl}&=
  q^{(\lambda,\lambda)}
    \sum_{p,t\in I}\ra^{jk}_{pt} z_{ip}z_{tl}\label{eq-projl1},\\
  \sum_{m,n,p,t\in I}\rc^{kl}_{mt}\ra^{jm}_{np}z_{in}z_{pt}&=
  q^{(\lambda,\lambda)}
    \sum_{p,t\in I}\ra^{jk}_{pt} z_{ip}z_{tl}\label{eq-projr1},\\  
  q^{(\lambda,\lambda)}\sum_{i,j\in I}C_{ij}z_{ij}&=1, \qquad \mbox{ where }
         C_{kl}:=\sum_{i\in I}(\rg^-)_{kl}^{ii}.\label{eq-spur}
\end{align}
Instead of (\ref{eq-projl1}) and (\ref{eq-projr1}) one can also write
\begin{align}
\sum_{m,n,p,t\in I}\rhm {}^{ij}_{nm}\ra^{mk}_{pt}z_{np}z_{tl}&=
q^{-(\lambda ,\lambda )}\sum _{p,t\in I}\ra {}^{jk}_{pt}z_{ip}z_{tl}
\label{eq-projl2},\\
\sum _{m,n,p,t\in I}\rcm {}^{kl}_{mt}\ra {}^{jm}_{np}z_{in}z_{pt}&=
q^{-(\lambda ,\lambda )}\sum _{p,t\in I}\ra {}^{jk}_{pt} z_{ip}z_{tl}
\label{eq-projr2}.  
\end{align}
As the right sides of (\ref{eq-projl2}), (\ref{eq-projr2}) coincide one
obtains
\begin{align}
\label{eq-refleq}
\sum_{m,n,p,t\in I}\rhm {}^{ij}_{nm}\ra^{mk}_{pt}z_{np}z_{tl}=
\sum _{m,n,p,t\in I}\rcm {}^{kl}_{mt}\ra {}^{jm}_{np}z_{in}z_{pt}.
\end{align}

As in the classical case
\begin{align}\label{Aqlim}
  \A^q_\lambda\cong \underrightarrow{\lim} V(n\lambda)^\ast\ot V(n\lambda).
\end{align}
Indeed, relation (\ref{c-lamclam})
allows one to write any $n$-fold product of the generators $z_{ij}$ as a linear
combination of $c_{F(v_\lambda)^{\ot n}}^{n\lambda}
  c_{G(f_{-\lambda})^{\ot n}}^{-nw_0\lambda}$, with
  $F\in V(n\lambda)^\ast$ and $G\in V(n\lambda)$. Therefore the arguments
  preceding (\ref{Alim}) can be repeated to verify (\ref{Aqlim}).
As $S_q[G/P_S]$ and $S_q[G/P_S^{op}]$ are quadratic algebras and
(\ref{Aqlim}) holds, Equations (\ref{eq-projl1}) -- (\ref{eq-spur}) form a
complete set of defining relations for $\A_\lambda^q$.

The $q$-deformed analogue of $\C[G/L_S]$  is defined by
\begin{align}\label{cqgldef}
  \C_q[G/L_S]=\{a\in \C_q[G]\,|\, a_{(1)}\,a_{(2)}(k)=\vep(k)a\quad\forall
  k\in K\},
\end{align}  
where $K:=\U_q(\lfrak_S)$ is the Hopf subalgebra of $\U$ generated by
the elements $\{K_i,K_i^{-1},E_j,F_j\,|\,i=1,\dots,r,\, j\in S\}$.
By construction $\C_q[G/L_S]$ is a left $\C_q[G]$-comodule algebra containing
$\A_\lambda^q$. The following proposition was proved for special cases
by M.~S.~Dijkhuizen and J.~V.~Stokman \cite{a-DS992} and in full generality by
J.~V.~Stokman \cite{a-stok02p} using $C^\ast$-algebra techniques.
Here we give an alternative proof based on the corresponding classical
result.

\begin{proposition}
  $\A_\lambda^q\cong \C_q[G/L_S]$ as left $\C_q[G]$-comodule algebras.
\end{proposition}

\begin{proof}
  Recall that the spaces $\C_q[G]$ and $\C[G]$ are isomorphic cosemisimple
  coalgebras with finite dimensional isotypical components
  (cf.~(\ref{peterw}) and
  (\ref{qpeterw})). The decompositions of $\C_q[G/L_S]$ and of $\C[G/L_S]$ into
  irreducible $\C[G]$-subcomodules coincide. This follows from the fact, that
  the dimensions of the weight spaces of irreducible representations $V(\mu)$
  of $\U$ are the same as in the classical case.
  On the other hand by (\ref{Alim}) and (\ref{Aqlim}) the $\C[G]$-comodules
  $\A_\lambda$ and $\A_\lambda^q$ are isomorphic. Now the assertion follows
  from Proposition \ref{A=CGL} and from $\A_\lambda^q\subset \C_q[G/L_S]$. 
\end{proof}

For $q\in \R$ it is well known that $\cqgl$ can be endowed with a 
$\ast$-structure induced by the compact real form of $\U$.
For general $q\in \C$ there remains a $\C$-linear algebra antiautomorphism,
which will prove useful in later arguments.

Consider the $\C$-linear algebra antiautomorphism
coalgebra homomorphism $\varphi_U:\U\rightarrow \U$ given by
\begin{align*}
  K_i\mapsto K_i,\, E_i\mapsto K_iF_i,
  \, F_i\mapsto E_i K_i^{-1}, \qquad i=1,\dots,r.
\end{align*}  
The composition
\begin{align*}
  \varphi:=\kappa\circ\varphi_U^\ast:\cqg\rightarrow\cqg,\quad
  a\mapsto \big(u\mapsto a(\varphi_U(\kappa(u)))\,\,\forall u\in U\big)
\end{align*}
is a $\C$-linear algebra antiautomorphism coalgebra homomorphism satisfying
$\varphi^2=\id$. As $K$ is a Hopf subalgebra of $U$ satisfying
$\varphi_U(K)\subset K$ one obtains $\varphi(\B)\subset \B$.
For any multi-index $J=(j_1,\dots,j_n)$ of length $|J|:=n$ of nonnegative
integers $j_i\le r$ define $E^J:=E_{j_1}\dots E_{j_n}$ and
$\wght(J):=\sum_{i=1}^n \alpha_{j_i}$,
similarly define $K^J$ and $F^J$. Note that $\varphi_U\circ \kappa$ is an
algebra homomorphism such that
\begin{align*}
  \varphi_U\circ \kappa(E^J)=(-1)^{|J|}F^J,\quad
  \varphi_U\circ \kappa(F^J)=(-1)^{|J|}E^J,\quad
  \varphi_U\circ \kappa(K^J)=(K^J)^{-1}.
\end{align*}  
Then for any irreducible finite dimensional
representation $V(\mu)$ with highest weight vector $v_\mu$ and any
$f\in V(\mu)^*$ one obtains
\begin{align*}
  \varphi(c^\mu_{f,v_\mu})(E^JK^{J'}F^{J''})=\delta_{J'',0}(-1)^{|J|}
  q^{-(\wght(J'),\mu)}f(F^Jv_\mu).
\end{align*}
On the other hand for the dual representation $V(\mu)^\ast\cong
V(-w_0\mu)$ with lowest weight vector $f_{-\mu}$ and any 
$v\in V(\mu)$
\begin{align*}
  c^{-w_0\mu}_{v,f_{-\mu}}(E^JK^{J'}F^{J''})=\delta_{J'',0}
  q^{-(\wght(J'),\mu)}v(E^Jf_{-\mu}).
\end{align*}
Recall that for $\mu=\sum_{i=1}^r n_i\omega_i$ the left ideal in $\Uqnm$ which
annihilates $v_\mu$ is generated by $\{F_i^{n_i+1}\,|\,i=1,\dots,r\}$
and the left ideal in $\Uqnp$ which annihilates $f_{-\mu}$ is generated
by $\{E_i^{n_i+1}\,|\,i=1,\dots,r\}$, \cite[4.3.6]{b-Joseph}.
Therefore there is a well defined isomorphism of vector spaces
$\tilde{\varphi}:V(\mu)^\ast\rightarrow V(\mu)$ such that 
\begin{align*}
  (E^J f_{-\mu})(\tilde{\varphi}(f))=(-1)^{|J|}f(F^Jv_\mu)
\end{align*}
for $f\in V(\mu)^\ast$ and all multi-indices $J$. Note that $\tilde{\varphi}$ maps weight spaces
to their duals and by definition
\begin{align*}
  \varphi(c^\mu_{f,v_\mu})&=
  c^{-w_0\mu}_{\tilde{\varphi}(f),f_{-\mu}}, &
  \varphi(c^{-w_0\mu}_{v,f_{-\mu}})&=
  c^\mu_{\tilde{\varphi}^{-1}(v),v_\mu}.
\end{align*}  
For later reference we collect the above considerations in the following lemma.
\begin{lemma}\label{antiauto}
  For $0\neq q\in\C$ not a root of unity there exist a $\C$-linear
  map $\varphi:\cqgl\rightarrow\cqgl$ and a vector space isomorphism
  $\tilde{\varphi}:V(\lambda)^\ast\rightarrow V(\lambda)$
  with the following properties:
  \begin{enumerate}
    \item The map $\varphi$ is an algebra antiautomorphism and a coalgebra
          automorphism.
    \item The map $\tilde{\varphi}$ maps weight spaces to their duals.
    \item The relation
             $\varphi(z_{fv})=z_{\tilde{\varphi}^{-1}(v)\tilde{\varphi}(f)}$
          holds for the generators $z_{fv}$ of $\cqgl$.
  \end{enumerate}        
\end{lemma}

\section{Finite dimensional irreducible graded representations of
         $\C_q[G/L_S]$}\label{findimirrgrad}
To shorten notation we will write $\B:=\C_q[G/L_S]$ from now on.
The right action of $\U $ endows $\B $ with a natural $\rlat$-grading
given by
\begin{align*}
\deg \,z_{ij}&=\wght (f_i)+\wght (v_j)=\wght(v_j)-\wght(v_i)
\in\rlat .
\end{align*}

Let $\B _+,\B _0$, and $\B _-$ denote the unital subalgebras of $\B $
generated by
the sets $\{z_{ij}\}$ with $\deg \,z_{ij}\lordne 0$,
$\deg \,z_{ij}=0$ and $0\lordne \deg \,z_{ij}$, respectively.
We will also write $\B _+^+:=\{b\in \B _+\,|\,\vep(b)=0\}$.

\begin{lemma}\label{l-Bdreieck}
$\B =\B _-\B _0\B _+$.
\end{lemma}

\begin{proof}
It suffices to find a filtration $\Fil $ on $\B $ and
$\lambda _{ijkl}\in \C\setminus\{0\}$ 
such that 
$z_{ij}z_{kl}=\lambda _{ijkl}z_{kl}z_{ij}$ holds
in the associated graded algebra.

Recall the total ordering $\lord$ on $\wlat$ defined at the beginning of
Section \ref{genflag}.
Consider the totally ordered abelian semigroup
\begin{align*}
  \cN=\{(k,\mu_1,\ldots ,\mu_k,\nu_1,\ldots ,\nu_k)\,|&\,k\in \N_0,
  \mu_i,\nu_i\in\wlat,\\
        & \mu_j\lordne \mu_i \mbox{ or }(\mu_j=\mu_i, \nu_j\lord \nu_i) \,
        \forall \,i<j\}
\end{align*}
with the lexicographic ordering with respect to the ordering $\lord$ of
$\wlat$. The sum of two elements of $\cN$ is defined by
\begin{align*}
  (k,\mu_1,\ldots ,\mu_k,\nu_1,\ldots ,\nu_k)+&(l,\mu_1',\ldots ,\mu_k',
  \nu_1',\ldots ,\nu_k')\\
  &=(k+l,\mu_1'',\ldots ,\mu_{k+l}'',\nu_1'',\ldots ,\nu_{k+l}'')
\end{align*}
where $(\mu_i'',\nu_i'')$, $i=1,\dots,k+l$, are the elements of
\begin{align*}
  \{(\mu_j,\nu_j),(\mu_p',\nu_p')\,|\, j=1,\dots,k \mbox{ and } p=1,\dots,l\}
\end{align*}  
in lexicographically increasing order.

There exists an $\cN$-filtration
$\Fil$ on $\B$ defined by (recall the remark at the end of the introduction)
\begin{align*}
  \deg_\Fil(z_{ij})=(1,\wght(v_i),\wght(v_j)).
\end{align*}  

By (\ref{eq-refleq}) one has
\begin{align*}
  z_{ij}z_{kl}=\sum_{\scriptstyle \parbox{6ex}{\scriptsize
      \centerline{$a,b,c,d,m,$} \centerline{$n,p,t\in I$}}}\rgm {}^{jk}_{ab}\rh ^{ia}_{mc}
  \rcm {}^{bl}_{dt}\ra ^{cd}_{np}z_{mn}z_{pt}.
\end{align*}
Suppose that $\wght(v_i)\lordne\wght(v_k)$ or $\wght(v_i)=\wght(v_k)$,
$\wght(v_j)\lordne\wght(v_l)$.
Since on the right hand side $\wght(v_k)\lord \wght(v_a)\lord\wght(v_m)$
holds in each nonzero summand by (\ref{rtabelle}),
one obtains $\deg _\Fil (z_{mn}z_{pt})<\deg _\Fil (z_{kl}z_{ij})$
whenever $m\ne a$ or $a\ne k$. Moreover if $m=a$ then $i=c$. Thus
again by  (\ref{rtabelle}) one has $\wght(v_i)=\wght(v_c)\lord\wght(v_p)$ and
therefore $\deg _\Fil (z_{mn}z_{pt})<\deg _\Fil (z_{kl}z_{ij})$ for $c\ne p$.
Finally, if $c=p$ then $n=d$ and therefore
$\wght(v_l)\lord\wght(v_d)=\wght(v_n)$. Thus
$\deg _\Fil (z_{mn}z_{pt})<\deg _\Fil (z_{kl}z_{ij})$ in all cases different
from $i=c=p$, $j=b=t$, $k=a=m$, $l=d=n$. This yields
\begin{align*}
z_{ij}z_{kl}=q^{(\wght (v_i)-\wght (v_j),\wght (v_k)+\wght (v_l))}z_{kl}z_{ij}
\end{align*}
up to terms of lower degree with respect to $\Fil $.
\end{proof}

Let $\B _\mu $, $\mu \in \rlat $, denote the subspace of $\B $
consisting of elements of $Q$-degree $\mu $.
In the following all irreducible finite dimensional graded representations $V$
of $\B $ will be determined. Here, a graded representation is a graded
vector space $V=\bigoplus _{\mu \in \rlat } V_\mu $ with a $\B$ action such
that $\B _\mu V_\nu \subset V_{\mu +\nu }$ for all $\mu ,\nu \in \rlat $.
A graded representation of $\B$ will be called an irreducible graded
representation if it does not possess any nontrivial invariant graded
subrepresentation.

\begin{lemma}\label{l-uhw}
Let $V$ be a finite dimensional irreducible graded representation of $\B$
and let $\lambda_0$ denote the highest weight of $V$ with respect to $\lordne$.
Then $\dim V_{\lambda_0}=1$. 
\end{lemma}

\begin{proof}
Let $V$ be a finite dimensional irreducible graded
representation with highest weight $\lambda _0\in Q$.
Then $\lambda _0+\deg \,z_{ij}\lordne \lambda _0$ for $\deg \,z_{ij}\lordne 0$
and hence $\B _+^+V_{\lambda _0}=0$. By Lemma \ref{l-Bdreieck} and the
irreducibility of $V$, $V_{\lambda _0}$ is an irreducible representation of
$\B _0$. Given $v\in V_{\lambda_0}$, $v\neq 0$, act with both sides of
Equation (\ref{eq-refleq}) with $\wght (v_k)=\wght (v_j)$,
$\wght (v_l)=\wght (v_i)$ on $v$.
Using $\wght (v_i)\lord \wght (v_m)\lord \wght (v_t)$ on the left
hand side, $\wght (v_j)\lord \wght (v_p)$ and $\wght (v_t)\lord \wght (v_k)
=\wght (v_j)$ on the right hand side, and $\B _+^+v=0$ one obtains the
relations
\begin{align}
z_{jk}z_{il}v=z_{il}z_{jk}v \qquad \text{whenever
$\wght (v_i)=\wght (v_l)$, $\wght (v_j)=\wght (v_k)$.}
\end{align}
Since the generators of $\B _0$ commute on $V_{\lambda _0}$ and $V_{\lambda _0}$
is an irreducible representation of $\B _0$, one has $\dim V_{\lambda _0}=1$.
\end{proof}

The following facts are needed in the sequel.

\begin{lemma}\label{l-weights}
Let $\mu ,\nu $ be weights of $V(\lambda )$.
\begin{itemize}
\item[(i)] $\mu \in W\lambda $ if and only if $(\mu ,\mu )=(\lambda ,\lambda )$.

\item[(ii)] If $\mu \in W\lambda $ then $(\mu ,\mu )=(\mu ,\nu )$ if and
  only if $\mu =\nu $.

\item[(iii)] 
  For any $\mu \in W\lambda $, $\mu \not= \lambda $, there exists
$i\in\{1,\ldots ,r\}$ with $(\mu ,\alpha _i)<0$. In this case  $\mu -\alpha _i$
is not a weight of $V(\lambda )$.
\end{itemize}
\end{lemma}

\begin{proof}
The first two statements can be found in \cite[11.4]{b-Kac1}

(iii) If $(\mu ,\alpha _i)\ge 0$ for all $i$ then $\mu $ is dominant.
As $\lambda$ is the only dominant weight in $W\lambda$ this proves the first
statement.
If $\mu=w\lambda$, $w\in W$, and $(\lambda, w^{-1}\alpha_i)=(\mu,\alpha_i)<0$
then $w^{-1}\alpha_i\in R^-$. Hence
$w^{-1}(\mu-\alpha_i)=\lambda-w^{-1}\alpha_i$ is not a weight of $V(\lambda)$
which proves the second statement. 
\end{proof}

Let $V^0=\C v$ denote the trivial representation, i.\,e.~$bv=\vep (b)v$
for all $b\in\B $.

\begin{proposition}\label{IrrGradRep}
Any finite dimensional irreducible graded representation of $\B $ is
isomorphic to $V^0$.
\end{proposition}

\begin{proof}
Let $V$ be such a representation and let $v\in V$ be a highest weight vector.
Lemma \ref{l-uhw} implies that for any 
$i,j\in I$ satisfying
$\wght (v_i)=\wght (v_j)$ there exists $\mu _{ij}\in \C $ such that
$z_{ij}v=\mu _{ij}v$. Moreover by (\ref{eq-spur}) there exist $i_0,j_0\in I$,
$\wght(v_{i_0})=\wght(v_{j_0})$, such
that $\mu _{i_0j_0}\not=0$. The proof is now performed in several steps.
The first term on the right hand side of (\ref{Rstruk}) and the properties
of the $R$-matrices collected in (\ref{rtabelle}) are frequently used.

\textit{Step 1. If $\wght (v_l)\not=\wght (v_{i_0})$ or
$\wght (v_{i_0})\lordne \wght (v_k)$ then $z_{kl}v=0$.}

The action of (\ref{eq-projl2}) and (\ref{eq-projr2}) on $v$ yields
\begin{align*}
(\ref{eq-projl2}),i=i_0,l=j_0,\wght (v_{i_0})\lordne \wght (v_j):\quad &
q^{(\wght (v_{i_0}),\wght (v_k)-\wght (v_j))}\mu _{i_0j_0}z_{jk}v=0,\\
(\ref{eq-projr2}),j=i_0,k=j_0,\wght (v_l)\lordne \wght (v_{i_0}):\quad &
\mu _{i_0j_0}z_{il}v=0.
\end{align*}
This implies that
\begin{align}\label{eq-zvnull}
z_{kl}v=0 \text{ for $\wght (v_{i_0})\lordne \wght (v_k)$ or $\wght (v_l)\lordne
\wght (v_{i_0})$.}
\end{align}
It remains to prove that $z_{il}v=0$ whenever
$\wght (v_{i_0})\lordne \wght (v_l)$ and $\wght (v_{i})\lord \wght (v_{i_0})$.
We proceed by induction over $\wght(v_i)$ with respect to $\lord$.
Suppose that $i\in I$ such that
$\wght (v_i)\lord \wght (v_{i_0})$ and $z_{al}v=0$
for all $a$ with $\wght (v_i)\lordne \wght (v_a)$.
By (\ref{eq-zvnull}) this is fulfilled  if $\wght (v_i)=\wght (v_{i_0})$.
For $i,l$ as above and any $j,k$ such that $\wght (v_k)\lordne \wght (v_j)$,
the action of
(\ref{eq-refleq}) on $v$ yields
$q^{(\wght (v_i),\wght (v_k)-\wght (v_j))}z_{jk}z_{il}v=0$
by induction hypothesis and $\B^+_+v=0$.
Thus $z_{il}v$ is a weight vector of $V$ with $\B _+^+z_{il}v=0$.
Since $0\lordne \deg\,z_{il}$, by Lemma \ref{l-Bdreieck} the
representation $\B z_{il}v=\B _-\B _0z_{il}v$ is a graded subrepresentation
of $V$ which does not contain $v$. Since $V$ is irreducible, $z_{il}v=0$.

\textit{Step 2. \begin{minipage}[t]{11.5cm} $\wght (v_{i_0})\in W\lambda $.
   In particular $i_0=j_0$.\end{minipage}}

Insert $i=j=i_0$, $k=l=j_0$ in (\ref{eq-projl2}). By Step 1 one gets
\begin{align*}
q^{(\wght (v_{i_0}),\wght (v_{j_0})-\wght (v_{i_0}))}z_{i_0j_0}^2v=
q^{-(\lambda ,\lambda )+(\wght (v_{i_0}),\wght (v_{j_0}))}z_{i_0j_0}^2v.
\end{align*}
Since $z_{i_0j_0}v=\mu _{i_0j_0}v$ and $\mu _{i_0j_0}\not=0$, one has
$(\lambda ,\lambda )=(\wght (v_{i_0}),\wght (v_{i_0}))$. Then Lemma
\ref{l-weights}(i) applies.

\textit{Step 3. $z_{kl}V=0$ for any $k,l$ with
 $\wght (v_{i_0})\lordne \wght (v_k)$ or $\wght (v_{i_0})\lordne \wght (v_l)$.}

To shorten notation we will write $\nu_i:=\wght (v_i)$ in this step. 
Equation (\ref{eq-refleq}) yields
\begin{align*}
q^{(\wght (v_i),\wght (v_k)-\wght (v_j))}z_{jk}z_{il}&+\sum _{\scriptstyle
\parbox{4ex}{\scriptsize\centerline{$n,p,t\in I$}\centerline{$\nu_i\pordne \nu_t$}}}
\alpha _{npt}z_{np}z_{tl}=\\
&=q^{(\wght (v_l),\wght (v_j)-\wght (v_k))}z_{il}z_{jk}
+\sum _{\scriptstyle \parbox{4ex}{\scriptsize \centerline{$n,p,t\in I$}
\centerline{$\nu_l\pordne \nu_n$}}}\beta _{npt}z_{in}z_{pt}
\end{align*}
for all $i,j,k,l$ and some complex numbers $\alpha _{npt},\beta _{npt}$
(depending on $i,j,k,l$). By induction over $\wght(v_l)$ with respect to
$\pord$ this implies that
\begin{align}\label{eq-zmt-z}
z_{il}z_{jk}&=q^{(\wght (v_i)+\wght (v_l),\wght (v_k)-\wght (v_j))}z_{jk}z_{il}
+\sum _{\scriptstyle \parbox{7ex}{\scriptsize \centerline{$m,n,p,t\in I$}
\centerline{$\nu_i\pord \nu_p,\nu_l\pord \nu_t,(p,t)\not=(i,l)$}}}
\alpha _{mnpt}z_{mn}z_{pt}
\end{align}
where $\alpha _{mnpt}\in \C $ depend on $i,j,k,l$.
Since $z_{il}v=0$ for $\wght (v_{i_0})\lordne \wght (v_i)$ or
$\wght (v_{i_0})\lordne \wght (v_l)$ Equation (\ref{eq-zmt-z}) immediately
gives $z_{il}V=z_{il}\B v=0$ in this case.

\textit{Step 4. If $\wght (v_{k_0})\pord\wght (v_{i_0})$ fails for some
  $k_0$ then $z_{k_0i_0}V=0$.}

Suppose that $\wght (v_{k_0})\pord\wght (v_{i_0})$ does not hold. Then the
same is true for any $k_1$ with $\wght (v_{k_0})\pordne \wght (v_{k_1})$.
Thus one can perform the proof by induction on $\wght(v_{k_0})$ with respect
to $\pord$. First let both sides
of (\ref{eq-projl2}) with $(i,j,k,l)=(k_0,i_0,i_0,i_0)$ act on $v$.
By $\wght (v_i)\pord\wght (v_t)$ and the induction hypothesis
on the left hand side and by $\wght (v_j)\pord\wght (v_t)$ and Steps 1
and 2 on the right hand side one gets
$z_{i_0i_0}z_{k_0i_0}v=z_{k_0i_0}z_{i_0i_0}v$.
Next let both sides of (\ref{eq-projl2}) with $(i,j,k,l)=(i_0,k_0,i_0,i_0)$
act on $v$. Now $\wght (v_i)\pord\wght (v_t)$ and Step 1
on the left and $\wght (v_k)\pord\wght (v_p)$ and Step 3 on the right side
imply that
$q^{(\wght (v_{i_0}),\wght (v_{i_0})-\wght (v_{k_0}))}z_{k_0i_0}z_{i_0i_0}v=
q^{-(\lambda ,\lambda )+(\wght (v_{k_0}),\wght (v_{i_0}))}z_{i_0i_0}z_{k_0i_0}v$.
These two equations give
$$ 0=(1-q^{2(\wght (v_{i_0}),\wght (v_{i_0})-\wght (v_{k_0}))})z_{k_0i_0}
z_{i_0i_0}v.
$$
Since $\mu_{i_0i_0}\neq 0$, Lemma \ref{l-weights}(ii) yields $z_{k_0i_0}v=0$.
Finally, set $i=k_0$ and $l=i_0$ in (\ref{eq-zmt-z}). Then from Step 3
together with the induction hypothesis one obtains $z_{k_0i_0}V=0$.

\textit{Step 5. One has $\wght(v_{i_0})=\lambda$.
  By Step 1 and (\ref{eq-spur}) this implies that
  $V{=}\C v{=}V^0$.}

Suppose that $\wght(v_{i_0})\neq \lambda$. By Lemma \ref{l-weights}(iii)
there exists $j\in \{1,2,\ldots ,r\}$ such that
$(\wght (v_{i_0}),\alpha _j)<0$ and therefore $E_j v_{i_0}\neq 0$.
In particular the second statement of Lemma \ref{l-weights}(iii) implies
that
$F_jv_{i_0}=0$. Since by Step 2 the weight space of $V(\lambda)$ of weight
$\wght(v_{i_0})$ is one-dimensional there exists a weight basis of
$V(\lambda)$ respecting the decomposition $V(\lambda)=\bigoplus_m V(\lambda)_m$
into irreducible $U_j$-modules. As the previous steps hold for an arbitrary
weight basis, we can assume that the conditions of Lemma \ref{l-divers} are
satisfied and that there exists $i_1\in I$ such that $v_{i_1}=E_j v_{i_0}$.

For $m\ge 0$ set $v_{(m)}:=z_{i_1i_0}^mv$.
Then by (\ref{eq-zmt-z}) with $i=l=k=i_0$, $j=i_1$ and Step 3 one obtains
for all $m\ge 0$ the equation
\begin{align}\label{eq-ziivk}
&z_{i_0i_0}v_{(m)}=q^{-2m(\wght (v_{i_0}),\alpha_j)}\mu _{i_0i_0}v_{(m)}.
\end{align}
Now let both sides of (\ref{eq-projr2}) with $(i,j,k,l)=(i_1,i_1,i_1,i_0)$
act on $v_{(m)}$. By Steps 3 and 4 and Lemma \ref{l-divers}(iii) one gets
\begin{align}
z_{i_1i_0}z_{i_1i_1}v_{(m)}=
q^{(\wght (v_{i_1}),\wght (v_{i_1}))-(\lambda ,\lambda )}
z_{i_1i_1}z_{i_1i_0}v_{(m)}
-(q^{2(\wght (v_{i_0}),\alpha_j)}{-}1)z_{i_1i_0}z_{i_0i_0}v_{(m)}.
\end{align}
Since $z_{i_1i_1}v=0$ by Step 1, one obtains from (\ref{eq-ziivk})
by induction over $m$ the formula
\begin{align*}
  z_{i_1i_1}v_{(m)}=q^{-2m(\wght (v_{i_0}),\alpha _j)}
 (q^{2(\wght (v_{i_0}),\alpha _j)}-1)(1-q^{-2md_j})/(q^{2d_j}-1)\mu _{i_0i_0}
  v_{(m)}.
\end{align*}

Consider (\ref{eq-refleq}) with $i=k=i_1$, $j=l=i_0$. Then in view of
the Steps 3 and 4 and Lemma \ref{l-divers}(iii) one obtains
\begin{align*}
q^{(\wght (v_{i_1}),\alpha _j)}z_{i_0i_1}z_{i_1i_0}v_{(m)}
+(q^{(\wght (v_{i_0}),\alpha _j)}{-}q^{-(\wght (v_{i_0}),\alpha _j)})
(z_{i_1i_1}z_{i_0i_0}v_{(m)}{-}z_{i_0i_0}^2v_{(m)})=&\\
=q^{-(\wght (v_{i_0}),\alpha _j)}z_{i_1i_0}z_{i_0i_1}v_{(m)}.
&
\end{align*}
Since $z_{i_0i_1}v=0$ by Step 1, one obtains by induction for all $m\ge 1$
the formula
\begin{align}\label{eq-z01-vk}
z_{i_0i_1}v_{(m)}=q^{-4(m-1)(\wght (v_{i_0}),\alpha _j)}(1-q^{-2(\wght (v_{i_0}),
\alpha _j)})\frac{1-q^{-2md_j}}{q^{2d_j}-1}\mu _{i_0i_0}^2 v_{(m-1)}.
\end{align}
Since $V$ is finite dimensional and $0\pordne \deg z_{i_1i_0}$ there exists
$m>0$ such that $v_{(m)}=0$ but $v_{(m-1)}\ne 0$. Then $z_{i_0i_1}v_{(m)}=0$
which is a contradiction to (\ref{eq-z01-vk}).
\end{proof}
       
\section{ The coalgebra $\overline{U}=U/UK^+$ }\label{U/UK+}

From now on for the remaining sections of this paper we restrict to the case
of irreducible flag manifolds.
Let $s\in\{1,\dots,r\}$ denote the missing index,
i.~e.~$S\cup\{\alpha_s\}=\pi$ and therefore $\lambda=\omega_s$.
Fix a reduced decomposition of the longest element $w_0\in W$ of the Weyl
group. Let $E_\beta, F_\beta$, $\beta\in R^+$, denote the corresponding root
vectors in $U$ \cite[8.1]{b-CP94}, \cite[6.2]{b-KS}.

\begin{lemma}\label{baslem}
 Let $g_i,i=1,2,\dots,\# R{+}r$, denote the generators $E_\beta,F_\beta,K_j$
 of $U$ with respect to $w_0$ in an arbitrary order. Then the elements
 \begin{equation*}
   \prod_{i=1}^{\# R{+}r} g_i^{n_i},
 \end{equation*}  
 $n_i\in \N_0$ if $g_i=E_\beta,F_\beta$ and $n_i\in\Z$ if $g_i=K_j$ form
 a vector space basis of $U$.
\end{lemma}
\begin{proof}
  This follows from the $q$-commutativity of the generators in the graded
  algebra associated to a certain filtration of $U$
  \cite[Prop.~1.7 d]{a-DCoKa90}.
\end{proof}  
Recall the decomposition $R^+=R_S^+\cup \gproots$ of the positive roots of
$\gfrak$. In the following we will use the abbreviation
$M:=\dim \gfrak/\pfrak_S=\#\overline{R^+_S}$.

\begin{proposition}\label{basprop}
  Let $\beta_1, \beta_2,\dots,\beta_{M}$ and
  $\beta'_1, \beta'_2,\dots,\beta'_{M}$ denote the elements
  of $\gproots$ in
  arbitrary fixed orders. Then the elements
  \begin{align}\label{basisUquer}
    \prod_{i=1}^M (E_{\beta_i})^{m_i} \prod_{j=1}^M (F_{\beta_j'})^{n_j},
  \end{align}
  $m_i,n_j\in\N_0$, form a vector space basis of $\Ubar=U/UK^+$.
\end{proposition}
\begin{proof}
  Let  $\gamma_1, \gamma_2,\dots,\gamma_{\# R^+_S}$ and
  $\gamma'_1, \gamma'_2,\dots,\gamma'_{\# R^+_S}$
  denote the elements of $R^+_S$ in
  arbitrary fixed orders. By Lemma \ref{baslem} the elements
  \begin{equation}
   \prod_{k=1}^M (E_{\beta_k})^{m_k}\prod_{l=1}^M (F_{\beta_l'})^{n_l}
   \prod_{i=1}^{\# R_S^+} (E_{\gamma_i})^{r_i} \prod_{j=1}^{\# R_S^+}
   (F_{\gamma_j'})^{s_j}
    K_1^{i_1}\dots K_r^{i_r}\label{basis},
  \end{equation}
  $i_k\in\Z,r_i,s_j,n_k,m_l\in\N_0$, form a vector space basis of $U$.
  Thus it suffices to show that the elements
  \begin{eqnarray}
     \prod_{k=1}^M (E_{\beta_k})^{m_k}\prod_{l=1}^M (F_{\beta_l'})^{n_l}
     (K_1^{i_1}\dots K_r^{i_r}-1),\label{bas1}\\
    \prod_{k=1}^M (E_{\beta_k})^{m_k} \prod_{l=1}^M (F_{\beta_l'})^{n_l}
    \prod_{i=1}^{\# R_S^+} (E_{\gamma_i})^{r_i}\prod_{j=1}^{\# R_S^+}
    (F_{\gamma_j'})^{s_j}
    K_1^{i_1}\dots K_r^{i_r},\label{bas2}
  \end{eqnarray}
   $i_k\in\Z,r_i,s_j,n_k,m_l\in\N_0$, $\sum_{i=1}^{\# R_S^+} (r_i+s_i)\ge 1$,
   form a vector space basis of $UK^+$.
   The expressions (\ref{bas1}) and (\ref{bas2}) form a set of
   linearly independent elements of $UK^+$. Any element of $UK^+$ can be
   written as a sum of expressions of the form
    \begin{eqnarray}
    \prod_{k=1}^M (E_{\beta_k})^{m_k} \prod_{l=1}^M (F_{\beta_l'})^{n_l}
    \prod_{i=1}^{\# R_S^+} (E_{\gamma_i})^{r_i}\prod_{j=1}^{\# R_S^+}
    (F_{\gamma_j'})^{s_j}
    K_1^{i_1}\dots K_r^{i_r}G\label{Gexpression}
  \end{eqnarray}
  where $G\in\{K_j-1,F_i,E_i|i\neq s\}$ and
  $i_k\in\Z,r_i,s_j,n_k,m_l\in\N_0$.
  For $G=E_i$, $i\neq s$, we show that (\ref{Gexpression}) is a linear
  combination of elements of the form (\ref{bas1}) and (\ref{bas2}).
  The cases $G=F_i$, $i\neq s$, and $G=K_i-1$ are dealt with in a similar way.

  First observe that $G=E_i$ $q$-commutes with $K_1^{i_1}\dots K_r^{i_r}$.
  If $\prod_{j=1}^{\# R_S^+}(F_{\gamma_j'})^{s_j}\neq F_i$ then by
  reordering $\prod_{j=1}^{\# R_S^+}(F_{\gamma_j'})^{s_j}E_i$ according to
  the above basis
  (\ref{basis}) one obtains a linear combination of monomials of the form
  $$ E_i^\delta \prod_{j=1}^{\# R_S^+}(F_{\gamma_j'})^{s_j'}
     K_1^{i_1'}\dots K_{r}^{i_r'}$$
  where $\delta=1$ or ($\delta=0$ and  $\sum_{j=1}^{\# R_S^+} s_j'> 0$).
  As the elements $E_{\gamma_i}$
  generate a subalgebra of $U$ with basis  $\prod_{i=1}^{\# R_S^+}
  (E_{\gamma_i})^{r_i}$ and $E_i$ is an element of this subalgebra the
  expression (\ref{Gexpression})
  for $G=E_i$ can indeed be written as a linear combination of elements
  of the form (\ref{bas2}). If on the other hand
  $\prod_{j=1}^{\# R_S^+}(F_{\gamma_j'})^{s_j}= F_i$ then the relation
  $$F_i E_i=E_i F_i-\frac{K_i-K_i^{-1}}{q_i-q_i^{-1}}$$
  implies the claim. 
\end{proof}

The canonical projection $\Phi:\U\rightarrow \Ubar$ is a surjective
coalgebra map. Recall that the coradical $U_0$ of $U$ is the subalgebra
generated by the elements $K_i,K_i^{-1},i=1,\dots,r$,
\cite[Lem.~5.5.5]{b-Montg93}.
For any surjective coalgebra map $f:C\rightarrow D$ the coradical of $D$ is
contained in the image of the coradical of $C$ \cite[Cor.~5.3.5]{b-Montg93}.
Therefore $\Ubar$ is connected, i.~e.~the coradical of
$\Ubar$ is one-dimensional. 

Set $\Ubar _+=\Lin_\C \{\prod _{i=1}^M E_{\beta _i}^{m_i}\}\subset \Ubar$,
$\Ubar _-=\Lin_\C \{\prod _{i=1}^M F_{\beta '_i}^{n_i}\}\subset \Ubar$,
where the products are
taken over all $i$ such that $\beta_i,\beta_i'\in\gproots$.
As  $\Ubar _+$ and $\Ubar _-$ are the images of the Hopf subalgebras
$\Uqbp$ and $\Uqbm$, respectively, under the canonical
projection $\Phi$, the subspaces $\Ubar _+\subset \Ubar$ and
$\Ubar _-\subset \Ubar$ are subcoalgebras. Note that $\Ubar_+$ and
$U/U(K^+ +\C F_s)$ are isomorphic coalgebras and hence $\Ubar_+$ admits a left
$U$-module structure. Similarly, the coalgebras $\Ubar_-$ and
$U/U(K^+ +\C E_s)$ are isomorphic.

Let $\Fil^s$ denote the filtration on $U$ defined by
\begin{align*}
  \deg^s(E_i)=\delta_{i,s}=\deg ^s(F_i),\quad \deg^s(K_i)=0.
\end{align*}
Then $\Fil^s$ induces a filtration on the  $U$-module $\Ubar$ by
\begin{align*}
  \deg ^s\left( \prod_{i=1}^M (E_{\beta_i})^{m_i}\prod_{j=1}^M
    (F_{\beta_j'})^{n_j}\right)
  =\sum_{i=1}^M m_i+n_i
\end{align*}
and on $\Ubar_+$ and $\Ubar_-$ all of which will
also be denoted by $\Fil^s$.
Further let $C$ denote the respective coradical filtrations.
Note that $\Fil^s_k\Ubar\subset C_k\Ubar$ and
$\Fil^s_k\Ubar_\pm\subset C_k\Ubar_\pm$.
On $\Ubar_\pm$ the filtration $\Fil^s$ is induced by a grading. The
homogeneous components of degree $k$ of this grading will be denoted by
$\Ubar_{\pm,k}$, i.\,e.~ $\Fil^s_k(\Ubar_\pm)=\bigoplus_{i=0}^k \Ubar_{\pm,i}$.

\begin{proposition}\label{coradpm}
  The coradical filtrations on $\Ubar_+$ and $\Ubar_-$ coincide with $\Fil^s$.
\end{proposition}
\begin{proof}
  We will prove the proposition for $\Ubar_+$. Define a grading $\cG$ on
  $\Uqbp$ by
  \begin{align*}
     \deg_\cG(E_i)=1,\quad \deg_\cG(K_i)=\deg_\cG(K_i^{-1})=0,
                        \qquad i=1,\dots,r.
  \end{align*}
  The induced grading on $\Ubar_+$ will also be denoted by $\cG$.
  Define linear functionals $\e_{\beta_i}\in(\Ubar_+)^\ast$ by
  \begin{align*}
     \e_{\beta_i}(E_{\beta_j})=\delta_{ij},\quad \e_{\beta_i}(X)=0
  \end{align*}  
  for all $X\in \Ubar_{+,k}$ with $k\neq 1$.
  Assume that the elements of $\gproots$ are ordered by increasing height,
  i.~e.
  \begin{align*}
    \hght(\beta_i)>\hght(\beta_j)\Rightarrow i>j.
  \end{align*}
  It suffices to verify the following statement:

  ($\ast$)\quad\begin{minipage}[t]{12cm}\textit{The set of functionals
      \begin{align*}
        \left\{\prod_{i=1}^M (\e_{\beta_i})^{n_i}\,\Bigg|\, (n_1,\dots,n_M)\in\N_0,
    \sum_{i=1}^M n_i=k\right\}
      \end{align*}
      is a basis of the homogeneous component $(\Ubar_{+,k})^\ast$
    of the graded dual algebra
  $\Ubar_+^{\gr\ast}=\bigoplus_{k=0}^\infty (\Ubar_{+,k})^\ast$.}
    \end{minipage}

  \noindent Indeed, $\Fil^s_k\Ubar_+\subset C_k\Ubar_+$. On the other hand
  if $x\in \Fil^s_k\Ubar_+\setminus\Fil^s_{k-1}\Ubar_+$ then by ($\ast$)
  there exists $f=\prod_{i=1}^M (\e_{\beta_i})^{n_i}$, $\sum_{i=1}^M n_i=k$,
  such that $f(x)\neq 0$. As $f|_{C_{k-1}\Ubar_+}=0$ this implies that
  $x\notin C_{k-1}\Ubar_+$.

 To verify ($\ast$) note first that as the elements of $\gproots$ are
 ordered by increasing height one has $\e_{\beta_i}(XE_{\beta_t})=0$ for
 all $X\in \Uqbp$, $i<t$,
 and $\e_{\beta_t}(XE_{\beta_t})=0$ for all $X\in \sum_{j=1}^r E_j\Uqbp$.
 This implies that for $1\le t\le M$, $m_t\ge n_t$ the relation
  \begin{align*}
    \left(\prod_{i=1}^t (\e_{\beta_i})^{n_i}\right)&
    \left(\prod_{i=1}^t (E_{\beta_i})^{m_i}\right)=
    q^{\left(n_t\beta_t,\sum_{i=1}^tm_i\beta_i-n_t\beta_t\right)}\times \\
    &\times\left(\prod_{i=1}^{t-1} (\e_{\beta_i})^{n_i}\right)
    \left(\prod_{i=1}^{t-1} (E_{\beta_i})^{m_i}E_{\beta_t}^{m_t-n_t}\right)
    \cdot\e_{\beta_t}^{n_t}\left( E_{\beta_t}^{n_t}\right)
  \end{align*}
  holds. For the same reason
  \begin{align*}
    \left( \prod_{i=1}^t (\e_{\beta_i})^{n_i}  \right)
    \left(\prod_{i=1}^{t+1} (E_{\beta_i})^{m_i}\right)=0
  \end{align*}
  if $m_{t+1}\neq 0$. Therefore
  \begin{align*}
    \left( \prod_{i=1}^{M} (\e_{\beta_i})^{n_i}  \right)
    \left(\prod_{i=1}^{M} (E_{\beta_i})^{m_i}\right)=0
  \end{align*}
  whenever $(n_{M},\dots,n_2,n_1)<
            (m_{M},\dots,m_2,m_1)$ with respect to the
  lexicographic ordering.            
  To prove ($\ast$) it now suffices to show that for all
  $\beta\in \overline{R^+_S}$ and $n\in \N$
  \begin{align*}
     \e_{\beta}^n\left(E_{\beta}^n  \right)
     =\e_{\beta}^{\ot n}(\kow^{n-1}(E_\beta)^n)\neq 0.
  \end{align*}
  For $\beta=\sum_{i=1}^r n_i\alpha_i\in Q$ define
  $K_\beta=\prod_{i=1}^r K_i^{n_i}$.
   To evaluate the above expression one verifies by induction over $n$ that
   \begin{align*}
     {\kow}^{n-1}(E_\beta)^n=
     \left[
       \sum
      _{i=1}^n
       1 \ot\dots\ot 1\ot
                  \overset{i\atop\downarrow}{E_\beta}\ot
                  K_\beta\ot\dots\ot K_\beta
                \right]^n+
         \sum_i
         X_{i1}\ot\dots\ot X_{in}
  \end{align*}
  where $X_{ij}\in \Uqbp$ and for all $i$ there exists $j\in \{1,\dots,n\}$
  such that $X_{ij}\in\sum_{k=1}^r \Uqbp E_k E_\beta +
  \sum_{ k\neq s}\Uqbp E_k$
  holds. Since  $K_\beta E_\beta=q^{(\beta,\beta)}E_\beta K_\beta$ and
  $e_\beta$ vanishes on $\sum_{k=1}^r\Uqbp E_kE_\beta$ one obtains
  \begin{align*}
     \e_{\beta}^n\left(E_{\beta}^n  \right)
     &=\e_{\beta}^{\ot n}\left(\Phi^{\ot n}\left(
                \left[\sum_{i=1}^n 1\ot\dots\ot1\ot
                  \overset{i\atop\downarrow}{E_\beta}\ot
     K_\beta\ot\dots\ot K_\beta\right]^n\right)\right)\\
     &=\prod_{k=1}^n\frac{q^{k(\beta,\beta)}-1}{q^{(\beta,\beta)}-1}
  \end{align*}
  where $\beta\in\overline{R^+_S}$ and
  $\Phi:\Uqbp\rightarrow \Ubar_+$ denotes the canonical projection.
  
\end{proof}

\begin{lemma}\label{invK}
$KC_l\Ubar \subset C_l\Ubar $ for all $l\ge 0$.
\end{lemma}
\begin{proof}
Since $K$ is a coalgebra, for $x \in\Ubar $, $k\in K$ we get
\begin{align*}
  (\kow k)(x\ot 1+1\ot x)=kx\ot 1+1\ot kx \in \Ubar\ot \Ubar. 
\end{align*}
Thus if $x\in C_l\Ubar$ and $k\in K$ then by \cite[Lem.~5.3.2(2)]{b-Montg93}
\begin{align*}  
\kopr (kx)-kx\ot 1-1\ot kx\in (\kopr k)\sum_{i=1}^{l-1}C_i\Ubar\ot
C_{l-i}\Ubar
\end{align*}
which proves the statement by induction over $l$.
\end{proof}
\begin{proposition}\label{corad}
  The coradical filtration of $\Ubar$ coincides with $\Fil^s$.
\end{proposition}

\begin{proof}
Recall that $\Fil^s_k\Ubar\subset C_k\Ubar$ and
$\kow x-1\ot x-x\ot 1\in \sum_{i=1}^{k-1}
\Fil_i^s\Ubar\ot \Fil^s_{k-i}\Ubar$ for all $x\in \Fil^s_k\Ubar$.
Contrary to the assertion of the proposition assume that $\Fil^s_k\Ubar\neq
C_k\Ubar$ for some $k\in \N_0$ and that $k$ is minimal with this property. 
Choose
$u=\sum \lambda _{m_1\dots m_M n_1\dots n_M}
\prod_{i=1}^M (E_{\beta_i})^{m_i} \prod_{j=1}^M (F_{\beta_j'})^{n_j} \in
C_k\Ubar\setminus \Fil^s_k\Ubar $.
Note that $\Ubar$ possesses a $\Z^r$-grading induced by the standard
$\Z^r$-grading of $U$.
By Lemma \ref{invK} one may assume that $u$ is homogeneous with respect to
the $\Z^r$-grading.
Moreover, without loss of generality we can assume that
$\lambda _{m_1\dots m_M n_1\dots n_M}=0$ whenever $\sum_{i=1}^M m_i+n_i\le k$.
By Proposition \ref{coradpm} one has $u\notin \Ubar_+\oplus\Ubar_-$.

Let $S_u\subset \N^{M}_0$ denote the subset
defined by
\begin{equation*}
  S_u:=\{(m_1,\dots,m_M)\,|\,\exists (n_1,\dots,n_M)\text{ such that }
                            \lambda_{m_1\dots m_M n_1\dots n_M}\neq 0\}.
\end{equation*}
Choose a multi-index $(k_1,\dots,k_M)\in S_u$ such that
  $\prod_{i=1}^M (E_{\beta_i})^{k_i}$ is maximal among the
  $\prod_{i=1}^M (E_{\beta_i})^{m_i}$,
  $(m_1,\dots,m_M)\in S_u$, with respect to the
  grading $\cG$ defined in the proof of Proposition \ref{coradpm}.
  The assumption $u\notin \Ubar_-$ implies that
  $\prod_{i=1}^M (E_{\beta_i})^{k_i}\neq 1$. Pick $(k_1',\dots,k_M')$ such that
  $\lambda_{k_1\dots k_M k_1'\dots k_M'}\neq 0$.
  Observe that then $\prod_{i=1}^M (F_{\beta_i'})^{k_i'}\neq 1$.
  Indeed, otherwise
  $u\in \Ubar_+$, since $u$ is homogeneous and
  $\prod_{i=1}^M (E_{\beta_i})^{k_i}\neq 1$ is maximal with respect to the
  grading $\cG$.
  Write $\kopr u\in \Ubar\ot \Ubar$ with respect to the basis given in
  Proposition \ref{basprop} in the first tensor factor. The second
  tensor factor corresponding to $\prod_{i=1}^M (F_{\beta_i'})^{k_i'}$ is
  given by
  \begin{equation*}
    \sum_{(m_1,\dots,m_M)\in S_u} \lambda _{m_1\dots m_M k'_1\dots k'_M }
     \prod_{i=1}^M (E_{\beta_i})^{m_i}\neq 0
  \end{equation*}   
  as $u$ is homogeneous with respect to the $\Z^r$-grading.
  But this means that
  \begin{align*}
    \kow(u)-1\ot u-u\ot 1\notin
    \sum_{i=1}^{k-1}\Fil^s_i\Ubar\ot \Fil^s_{k-i}\Ubar=
    \sum_{i=1}^{k-1} C_i\Ubar\ot C_{k-i}\Ubar.
  \end{align*}
  This contradicts $u\in C_k\Ubar$.
\end{proof}  

\section{Duality}\label{duality}

By definition of $\B=\cqgl$ the nondegenerate pairing (\ref{qpair}) induces a
well defined pairing
\begin{align}\label{qBpair}
  \pair{\cdot}{\cdot}:\B\ot \Ubar\rightarrow \C
\end{align}
of the algebra $\B$ and the coalgebra $\Ubar$.
\begin{proposition}\label{BUnondeg}
  The pairing (\ref{qBpair}) is nondegenerate.
\end{proposition}

\begin{proof}
  Since (\ref{qpair}) is nondegenerate $\Ubar$ separates $\B$ in
  (\ref{qBpair}). To show that $\B$ separates $\Ubar$ take any
  $x\in \Ubar\setminus\{0\}$. By Proposition \ref{basprop} the element
  $x$ can be represented by linear combination of elements of $U$ of
  the form (\ref{basisUquer}). Choose a basis vector
  \begin{align*}
    \prod_{i=1}^M (E_{\beta_i})^{m_i} \prod_{j=1}^M (F_{\beta_j'})^{n_j}\in U
  \end{align*}
  occurring with non-vanishing coefficient $a_0$ in such a linear combination
  such that $(\sum_{j=1}^M n_j,\sum_{i=1}^M m_i)$
  is maximal with respect to the lexicographic order.
  For given $n\in \N$ the elements
  $\prod_{i=1}^M (E_{\beta_i})^{k_i}f_{-n\omega_s}$, $\sum_{i=1}^M k_i\le n$,
  where $f_{-n\omega_s}$ denotes a lowest weight vector of $V(n\omega_s)^\ast$,
  are linearly independent in $V(n\omega_s)^\ast$.
  Indeed, by \cite[4.3.6]{b-Joseph} one has 
  \begin{align*}
    V(n\omega_s)\cong
  \Uqnp\bigg/
  \left(\sum_{i\neq s} \Uqnp E_i +\Uqnp E_s^{n+1}\right)\cong \Ubar_+/
  \Uqnp E_s^{n+1}
  \end{align*}
  as $\Uqnp$-modules, where $\Uqnp E_s^{n+1}\subset \Ubar_+$.
  By similar reasoning the elements
  $\prod_{j=1}^M (F_{\beta_j'})^{l_j}v_{n\omega_s}$,
  $\sum_{j=1}^M l_j \le n$, where $v_{n\omega_s}$ denotes a highest weight
  vector of $V(n\omega_s)$, are linearly independent in $V(n\omega_s)$.
  
  \noindent Thus for any $n>\max\{\sum _{j=1}^M n_j,\sum_{i=1}^M m_i\}$
  one obtains
  \begin{align*}
    x(v_{n\omega_s}\ot f_{-n\omega_s})=
       a_0\left( \prod_{j=1}^M (F_{\beta_j'})^{n_j} v_{n\omega_s}\right)\ot
            \left(\prod_{i=1}^M (E_{\beta_i})^{m_i}f_{-n\omega_s} \right)+\dots
  \end{align*}
  where $\dots$ denote terms which are linearly independent of the first
  expression. As
  $c^{n\omega_s}_{f, v_{n\omega_s}}c^{-nw_0\omega_s}_{v, f_{-n\omega_s}}\in \B$
  for all $f\in V(n\omega_s)^\ast$ and $v\in V(n\omega_s)$ one obtains
  $\pair{\B}{x}\neq 0$.
\end{proof}
To shorten notation we introduce index sets
\begin{align*}
  I_{(k)}:=\left\{i\in I\,|\,(\omega_s-\wght(v_i),\omega_s)=k d_s
      \right\}.
\end{align*}
Moreover, assume that the basis vector $v_N$ is a highest weight vector of
$V(\omega_s)$. Recall that $M=\dim \gfrak/\pfrak_s=\# \gproots$.
\begin{lemma}\label{b+/b+2dim}For all $i,j$ such that
  $(i,j)\notin (I_{(1)}\times\{N\})\cup(\{N\}\times I_{(1)})$
  one has
  \begin{align*}
  z_{ij}-\delta_{iN}\delta_{jN}\in (B^+)^2.
  \end{align*}
  In particular $\dim \B^+/(\B^+)^2\le 2 M$. 
\end{lemma}
\begin{proof}
  By definition the generators $z_{ij}$ of $\B$ satisfy the relation
  $\sum_{k\in I} z_{ik}z_{kj}=z_{ij}$. For $i\neq N\neq j$ or $i=N=j$ using
  $\vep(z_{ij})=\delta_{iN}\delta_{jN}$ one obtains
  \begin{align*}
    z_{ij}-\delta_{iN}\delta_{jN}\in (\B^+)^2.
  \end{align*}
  Since
  $V(\omega_s)=\Lin_\C\{\prod_{i=1}^M(F_{\beta_i})^{n_i}v_N\,|\,\beta_i\in
  \gproots$, $n_i\in\N_0\}$ it suffices to show that
  $z_{Nj},z_{jN}\in (\B^+)^2$ if $\wght(v_j)=\omega_s-2\alpha_s-\beta$ for
  some $\beta\in Q$, $\beta \pord 0$.

  In view of (\ref{rtabelle}) relation (\ref{eq-projl1}) for $i,j<N$,
  $k=l=N$ implies that
  \begin{align}\label{kombinbp}
    \sum_{n\in I} \rh^{ij}_{nN}z_{nN}\in (\B^+)^2.
  \end{align}
  More precisely, 
  \begin{align}\label{kombinbp2} 
    \sum_{n\in I} \rh^{ij}_{nN}z_{nN}{+}\sum_{n,m\in I\setminus\{N\}}
      q^{(\omega_s,\wght(v_m)-\omega_s)}\rh^{ij}_{nm}z_{nN}z_{mN}-
      q^{(\omega_s,\wght(v_j))}z_{iN}z_{jN}\in (\B^+)^3.
  \end{align}  
  Consider now $v_n\in V(\omega_s)$ with
  $\wght(v_n)=\omega_s-2\alpha_s-\beta$ for some
  $\beta\in Q$, $\beta \pord 0$. We will show that there exist
  $i,j<N$ such that $\rh^{ij}_{nN}\neq 0$.
  This implies that
  \begin{align}\label{mapinjective}
   (P_<\ot P_<)\circ \rh_{\omega_s,\omega_s}:
     V_{\omega_s-2\alpha-\beta}\ot\C v_N\rightarrow P_< V(\omega_s)\ot
     P_<V(\omega_s)
  \end{align}
  is injective. Here $P_<: V(\omega_s)\rightarrow V(\omega_s)$ denotes the
  projection given by
  \begin{align*}
    P_<|_{V_\mu}:=(1-\delta_{\mu,\omega_s})\id|_{V_\mu}
  \end{align*}  
  for any weight space $V_\mu\subset V(\omega_s)$.
  Therefore $z_{nN}\in (\B^+)^2$.

  If $\rh^{ij}_{nN}= 0$ for all $i,j<N$ then
  \begin{align}\label{inVv}
     \rh_{\omega_s,\omega_s}(v_n\ot v_N)-q^{(\wght(v_n),\omega_s)}v_N\ot v_n
     \in V(\omega_s)\ot v_N.
  \end{align}
  As $\rh_{\omega_s,\omega_s}$ is a $\U$-module homomorphism $E_k v_n$
  also satisfies (\ref{inVv}) for all $k\neq s$. Therefore we can assume that
  $E_k v_n=0$ for all $k\neq s$. But in this case
  $E_sv_n\neq 0$ and with $F_s v_N\neq 0$ one obtains from (\ref{Rstruk})
  \begin{align*}
    \rh^{ij}_{nN}=q^{(\omega_s-\alpha_s,\wght(v_n)+\alpha_s)}
  (q^{d_s}-q^{-d_s})\neq 0 
  \end{align*}
  for $v_i=F_sv_N$ and $v_j=E_sv_n.$

  To verify $z_{Nn}\in (\B^+)^2$ for all $n$ with
  $\wght(v_n)=\omega_s-2\alpha_s-\beta$ and $\beta\in Q$, $\beta \pord 0$,
  apply the algebra antiautomorphism $\varphi$ from Lemma \ref{antiauto}
  to $z_{nN}$.
\end{proof}

The relations $V(\omega_s)=\Uqnm v_N$ and $V(\omega_s)^\ast=\Uqnp f_N$ induce
direct sum decompositions
\begin{align}
  V(\omega_s)=\bigoplus _k V(\omega_s)_{(k)},\qquad
  V(\omega_s)^\ast=\bigoplus _k V(\omega_s)^\ast_{(k)}
\end{align}
where
\begin{align*} V(\omega_s)_{(k)}&:=
       \Lin_\C\left\{\left(\prod_{j=1}^M (F_{\beta_j'})^{n_j}\right)v_N\,
         \Big|\,\sum_{j=1}^M n_j=k\right\}=\Lin_\C\{v_i\,|\,i\in I_{(k)}\}\\
               V(\omega_s)_{(k)}^\ast&:=
       \Lin_\C\left\{\left(\prod_{i=1}^M (E_{\beta_i})^{n_i}\right)f_N\,
         \Big|\,\sum_{i=1}^M m_i=k\right\}=\Lin_\C\{f_i\,|\,i\in I_{(k)}\}\\.
\end{align*}
and $\beta_i,\beta_j'$ denote the elements of $\gproots$ in the order fixed
in the beginning of Section \ref{U/UK+}. 

Note that (\ref{eq-projl1}) for $i=l=N$ and $j,k<N$ yields
\begin{align*}
  z_{jk}=q^{(\omega_s,\omega_s-\wght(v_j)-\wght(v_k))}
  \sum_{p,t\in I\setminus\{N\}}\ra^{jk}_{pt}z_{Np}z_{tN} \quad \mod \,(\B^+)^3.
\end{align*}
Using $\sum_{m\in I} z_{jm}z_{mk}=z_{jk}$ one obtains
\begin{align}\label{pNNt}
  z_{jN}z_{Nk}= q^{(\omega_s,\omega_s-\wght(v_j)-\wght(v_k))}
  \sum_{p,t\in I\setminus\{N\}} \ra^{jk}_{pt}z_{Np}z_{tN} \quad
  \mod \,(\B^+)^3.
\end{align}
Therefore Lemma \ref{b+/b+2dim} leads to an isomorphism of $\Z^r$-graded
vector spaces
\begin{align*}
  (\B^+)^2/(\B^+)^3=D^{2,0}\oplus D^{1,1}\oplus D^{0,2}
\end{align*}
where
\begin{align*}
  D^{2,0}&=\Lin_\C\{z_{iN}z_{jN}\in (\B^+)^2/(\B^+)^3\,|\,i,j\in I_{(1)}\}\\
  D^{1,1}&=\Lin_\C\{z_{Ni}z_{jN}\in (\B^+)^2/(\B^+)^3\,|\,i,j\in I_{(1)}\}\\
  D^{0,2}&=\Lin_\C\{z_{Ni}z_{Nj}\in (\B^+)^2/(\B^+)^3\,|\,i,j\in I_{(1)}\}.  
\end{align*}
Recall from Section \ref{U/UK+} that $\Ubar_{+,i}$ denotes the homogeneous
component of degree $i$ of the graded coalgebra $\Ubar_+$. 

\begin{lemma}\label{D20estimate}
  $\dim D^{2,0}\le \dim \Ubar_{+,2}$.
\end{lemma}  
\begin{proof}
  Consider the vector space
  \begin{align*}
     \Vtil:=&\Lin_\C \{ f_i{\ot} f_j \in V(\omega_s)^\ast{\ot} V(\omega_s)^\ast
    \,|\, (2\omega_s-\wght (v_i)-\wght(v_j),\omega_s)=2 d_s\}\\
    =&V(\omega_s)_{(1)}^\ast\ot V(\omega_s)_{(1)}^\ast\oplus
      V(\omega_s)_{(0)}^\ast\ot V(\omega_s)_{(2)}^\ast\oplus
      V(\omega_s)_{(2)}^\ast\ot V(\omega_s)_{(0)}^\ast
  \end{align*}
and the linear maps
\begin{align*}
 \Psi_1 & :\Vtil\rightarrow \Vtil,&  \Psi_1(f)&:= f\circ
 \rh_{\omega_s,\omega_s} -q^{(\omega_s,\omega_s)} f,\\
 \Psi_2 &: \Vtil \rightarrow D^{2,0},& \Psi_2(f_i\ot f_j)&:=
        \begin{cases}
          q^{(\omega_s,\wght (v_j))}z_{jN}&\mbox{ if } i=N, j\in I_{(2)},\\
          q^{(\omega_s,\omega_s)}z_{iN}&\mbox{ if } i\in I_{(2)},j=N,\\
          q^{(\omega_s,\wght (v_j))}z_{iN}z_{jN}&\mbox{ if } i,j\in I_{(1)}.
        \end{cases}  
\end{align*}
Recall from (\ref{kombinbp2}) and from the injectivity of
(\ref{mapinjective}) that $z_{jN}\in D^{2,0}$ for
all $j\in I_{(2)}$ and therefore $\Psi_2$ is well defined. Moreover,
$\Psi_2$ is surjective.
Now we claim that 
\begin{align}\label{psipsi}
  \Psi_2\circ\Psi_1=0.
\end{align}
Indeed, if $i,j\in I_{(1)}$ then  in $(\B^+)^2/(\B^+)^3$ one calculates using
(\ref{eq-projl1}) and Lemma \ref{b+/b+2dim}
\begin{align*}
  \Psi_2\circ\Psi_1(f_i\ot f_j)&=\Psi_2\left(\sum_{n,m\in I}\rh^{ij}_{nm}
    f_n\ot f_m
              -q^{(\omega_s,\omega_s)}f_i\ot f_j\right)\\
          &= \sum_{n\in I_{(2)}} \rh^{ij}_{nN}q^{(\omega_s,\omega_s)}  z_{nN}+
             \sum_{n,m\in I_{(1)}}
             q^{(\omega_s,\wght(v_m))}\rh^{ij}_{nm}z_{nN}z_{mN}\\
             &\qquad-q^{(\omega_s,\omega_s+\wght(v_j))}z_{iN} z_{jN}\\
          &= \sum_{n,m,p,t\in I} \rh^{ij}_{nm}\ra^{mN}_{pt}z_{np}z_{tN}
              -q^{(\omega_s,\omega_s)} \sum_{p,t\in I}
              \ra^{jN}_{pt}z_{ip} z_{tN}
             =0. 
\end{align*}  
The cases $i=N$, $j\in I_{(2)}$ and $i\in I_{(2)}$, $j=N$ are dealt with in a
similar manner.

Note that for $V(\mu)^\ast\subset V(\omega_s)^\ast\ot V(\omega_s)^\ast$
the restriction $\Psi_1|_{\Vtil\cap V(\mu)^\ast}$
is multiplication by a nonzero scalar if $\mu\neq 2\omega_s$ and $0$ if
$\mu=2\omega_s$. As $f_N\ot f_N$ is the lowest weight vector of the submodule
$V(2\omega_s)^\ast \subset V(\omega_s)^\ast \ot V(\omega_s)^\ast$ one obtains
$\Vtil\cap V(2\omega_s)^\ast=\Ubar_{+,2}(f_N\ot f_N)$, where $\Ubar_+$ is
interpreted as $\Uqnp/\sum_{i\neq s}\Uqnp E_i$.
Therefore $\dim\ker \Psi_1\le \dim \Ubar_{+,2}$. Combining this estimate
with (\ref{psipsi}) one gets
\begin{align*}
  \dim D^{2,0}&=\dim \im \Psi_2=\dim \Vtil-\dim\ker \Psi_2\\
  &\overset{(\ref{psipsi})}{\le} \dim \Vtil-\dim\im \Psi_1
  =\dim \ker \Psi_1 \le \dim\Ubar_{+,2}.
\end{align*}  
\end{proof}  

\begin{proposition}\label{BB+CkUpair}
  The pairing $\pair{\cdot}{\cdot}:\B/(\B^+)^{k+1}\ot C_k\Ubar\rightarrow \C$
  is nondegenerate.
\end{proposition}
\begin{proof}
  By Proposition \ref{BUnondeg} $\B$ separates $C_k\Ubar$. Therefore it
  suffices to verify that $\dim \B/(\B^+)^{k+1}\le\dim C_k\Ubar$. By
  (\ref{pNNt}) one obtains a decomposition
  \begin{align*}
    (\B^+)^k/(\B^+)^{k+1}=\sum_{i=0}^k D^{i,k-i}
  \end{align*}
  where $D^{i,k-i}=\Lin_\C\{\prod_{m=1}^{k-i}z_{Nj_m}\prod_{m=1}^i z_{l_mN}\,
         |\,j_m,l_m\in I_{(1)}\}.$
  Thus it suffices to show that
  \begin{align}\label{DiUi}
    \dim D^{0,i}\le\dim\Ubar_{+,i}.
  \end{align}
  The estimate $\dim D^{i,0}\le\dim\Ubar_{-,i}$ then follows by application
  of the antiautomorphism $\varphi$ from Lemma \ref{antiauto}.

  For $i=1$ equality holds in (\ref{DiUi}) as both $\Lin_\C\{z_{iN}\,|\, i\in
  I_{(1)}\}$ and $\Ubar_{-,1}$ are irreducible $K$-modules with highest weight
  $-\alpha_s$.
  To verify (\ref{DiUi}) for general $i$
  consider an ordering $\vdash$ of the set
  $I_{(1)}$ such that
  \begin{align*}
    \hght(\omega_s-\wght (v_i) )> \hght(\omega_s-\wght (v_j) )
              \,\Longrightarrow\, i\vdash j.
  \end{align*}
  Note then that the elements
  \begin{align}\label{D02basis}
    \{ z_{Ni}z_{Nj}\in (\B^+)^2/(\B^+)^3\,|\,i,j\in I_{(1)}, j\vdash i\}
  \end{align}
  form a basis of $D^{0,2}$. Indeed, recall the functionals
  $\e_{\beta}\in(\Ubar_+)^\ast$ defined in the proof of Proposition
  \ref{coradpm}, and note that $z_{Ni}|_{\Ubar_+}$ is a nonzero multiple of
  $\e_{\beta}$ if $E_\beta f_N=f_i$, $i\in I_{(1)}$. Hence by ($\ast$) from
  the proof of Proposition
  \ref{coradpm} the set (\ref{D02basis}) consists of $\dim \Ubar_{+,2}$
  linear independent elements which by Lemma \ref{D20estimate} span
  $D^{0,2}$. Thus the elements
  \begin{align*}
  \left\{\prod_{m=1}^i z_{Nj_m}\,\Big|\, j_m\in I_{(1)},\,j_{m+1}
    \vdash j_m\right\} 
  \end{align*}  
  span $D^{0,i}$ and (\ref{DiUi}) follows from Proposition \ref{basprop}.
\end{proof}  

In view of the results of the previous sections the proposition above
allows us to refine the duality between $\B$ and $\Ubar$ of Proposition
\ref{BUnondeg}. Let $\B^\circ$ denote the dual coalgebra of $\B$, i.e.~the
coalgebra generated by the matrix coefficients of all finite dimensional
representations of $\B$ (cf.~e.\,g.~\cite[Sect.~1.2]{b-Montg93}).
As $\B$ is a right $U$-module algebra the dual
coalgebra $\B^\circ$ is a left $U$-module coalgebra.
In the simplest case of quantized
$\C P^1$, the so called "Podle\'s' quantum sphere", the dual coalgebra
has been determined in \cite{a-heko02p}. The left $U$-action on
$\B^\circ$ permits to ask for elements of $\B^\circ$ with an additional
finiteness property, which we will call the \textit{locally finite part}
\begin{align*}
  F(\B^\circ,K):=\{f\in\B^\circ\,|\, \dim(Kf)<\infty\}.
\end{align*}
Note that $\Ubar\subset  F(\B^\circ,K)$ via the nondegenerate pairing
(\ref{qBpair}), cf.~Lemma \ref{invK}. The first statement of the following
theorem is the main result of this paper. Together with the second statement
which is merely a reformulation of the definition of $\B$ it furnishes
the duality between $\B$ and $\Ubar$. Recall that $\Ubar$ is a left
$U$-module coalgebra and therefore the dual algebra $\Ubar^\ast$ of linear
functionals on $\Ubar$ is a right $U$-module algebra.
\begin{theorem}\label{locfin}
  \begin{minipage}[t]{6cm} 1) $F(\B^\circ,K)=\Ubar$\\
                        2) $\B=\{b\in \Ubar^\ast\,|\,\dim(bU)<\infty\}$.
  \end{minipage}                      
\end{theorem}  

\begin{proof}
1) Let $K'\subset K$ denote the subalgebra generated by the elements
$K_i,K_i^{-1}$,
$i=1,\dots,r$. As $F(\B^\circ,K)\subset F(\B^\circ,K')$ it suffices to
show that $F(\B^\circ,K')\subset \Ubar$. Note that the coradical filtration
of $F(\B^\circ,K')$ is invariant under the left $K'$-action, i.~e.~
\begin{align*}
  K' C_n F(\B^\circ,K')\subset C_n F(\B^\circ,K'),
\end{align*}
as the generators $K_i, K_i^{-1}$ of $K'$ act on $F(\B^\circ,K')$ by coalgebra
automorphisms. Thus by Proposition \ref{IrrGradRep} the coalgebra
$F(\B^\circ,K')$ is connected,
i.\,e.\ $C_0 F(\B^\circ,K')=\C\vep$, and therefore
any $f\in C_n F(\B^\circ,K')$, $n\in \N_0$, vanishes on $(\B^+)^{n+1}$.
By Proposition \ref{BB+CkUpair} this implies that $f\in C_n\Ubar$.

2) The dual Hopf algebra $U^\circ$ of $U$ satisfies
\begin{align}\label{Ucirc}
  U^\circ=\{a\in U^\ast\,|\, \dim(aU)<\infty\}
\end{align}  
and contains $\cqg$ as the linear span of the matrix coefficients of the
representations $V(\mu)$, $\mu\in P^+$. Recall that $U$ is semisimple and
any irreducible representation of $U$ can be obtained by tensoring some
$V(\mu)$ with a one dimensional representation $D_\nu$, $\nu\in\{-1,1\}^r$,
given by  $K_i v=\nu_i v$ for $v\in D_\nu$. Therefore
\begin{align}\label{BinUcirc}
  \B=\{a\in U^\circ\,|\, a_{(1)}\,a_{(2)}(k)=\vep(k)a \quad\forall k\in K\}.
\end{align}
Inserting (\ref{Ucirc}) in (\ref{BinUcirc}) leads to the desired expression. 
\end{proof}

Duality of $\B$ and $\Ubar$ also holds in the graded setting. Note that the
algebra $\B$  admits a decreasing filtration given by $\Fil_n\B=(\B^+)^n$.
Let $\gr\B=\bigoplus_{n=0}^\infty (\B^+)^n/(\B^+)^{n+1}$ denote the associated
graded algebra. Define $(\gr \B)_+$ and $(\gr \B)_-$ to be the subalgebras
of $\gr \B$ generated by $D^{0,1}$ and $D^{1,0}$, respectively.
Let $(\gr \B)_{\pm,n}$ denote the elements of $(\gr \B)_\pm$ of degree $n$.

On the other hand let $\gr \Ubar$ denote the graded coalgebra
associated to the coradical filtration of $\Ubar$. Let further
$(\gr\Ubar)^{\gr\ast}$ denote the graded dual of $\gr \Ubar$, similarly
$\Ubar_+^{\gr\ast}$, $\Ubar_-^{\gr\ast}$.
Recall that the right $K$-module structure of $\B$ and the left $K$-module
structure of $U$ are compatible with the filtration $\Fil$ of $\B$ and the
coradical filtration of $\Ubar$, respectively.
The following Corollary is an immediate consequence of Proposition
\ref{BB+CkUpair}.
\begin{cor}
  The pairing (\ref{qBpair}) induces isomorphisms
  \begin{align*}
    \gr \B&\cong (\gr\Ubar)^{\gr\ast}&
    (\gr \B)_+&\cong \Ubar_+^{\gr\ast}&
    (\gr \B)_-&\cong \Ubar_-^{\gr\ast}
  \end{align*}
  of graded right $K$-module algebras.
\end{cor}  

The $K$-module algebras $(\gr \B)_+$ and $(\gr \B)_-$ have been constructed
in \cite{a-SinVa} as $q$-analogues of the polynomial algebra on the
prehomogeneous space $\gfrak_{-1}$. Note that the left $U$-module structures on
$\Ubar_+$ and $\Ubar_-$ induce right $U$-module structures on
$(\gr \B)_+$ and $(\gr \B)_-$, respectively.

By statement ($\ast$) from the proof of Proposition (\ref{coradpm}) the algebra
$(\gr \B)_+$ is generated by the elements $(\gr \B)_{+,1}$ of
degree one and quadratic relations. The $K$-module $(\gr \B)_{+,1}$ is
irreducible and  each weight space of $(\gr \B)_{+,1}$ is one-dimensional.
Therefore we have a decomposition
\begin{align}\label{grotgr}
  (\gr \B)_{+,1}\ot(\gr \B)_{+,1}=\bigoplus_i V_i
\end{align}
into irreducible $K$-modules $V_i$ and each of these modules occurs with
multiplicity one. The dimension of each weight space of $(\gr \B)_{+,2}$
coincides with the dimension of the corresponding weight space of the
symmetric elements of $V(\omega_s)_{(1)}\ot V(\omega_s)_{(1)}$.
Therefore one obtains a complete set of defining
relations for $(\gr \B)_+$ if one sets the "antisymmetric" components of
the tensor product (\ref{grotgr}) equal to zero.
More precisely, let $V(\omega_s)_{(1)}\ot V(\omega_s)_{(1)}=S_{(1)}\oplus
A_{(1)}$ denote the decomposition into "symmetric" and "antisymmetric"
$\U$-submodules. Then by the above arguments the following
statement is proved.
\begin{cor}
  The algebra $(\gr \B)_+$ is isomorphic to the quotient of the tensor
  algebra $\bigoplus_{k=0}^\infty V(\omega_s)_{(1)}^{\ot k}$ by the ideal
  generated by $A_{(1)}\subset V(\omega_s)_{(1)}^{\ot 2}$. Moreover,
  \begin{align*}
    \dim (\gr \B)_{+,k} ={ M+k-1\choose k}.
  \end{align*}  
\end{cor}  

An analogous result holds for $(\gr \B)_-$. In particular one obtains the
following statements
(cp.~the list \cite[p.~27]{b-BastonEastwood}, simple roots are ordered as in
\cite[p.~58]{b-Humphreys}).
\begin{enumerate}
\item If $\gfrak=\slfrak_N$ and $S=\pi\setminus\{\alpha_s\}$ then
  $(\gr\B)_-$ is isomorphic to the $U_q(\slfrak_s\times
  \slfrak_{N-s})$-module algebra $\C_q[\mbox{Mat}_{N-s,s}]$ of quantized
  $(N-s,s)$-matrices.
  Similarly $(\gr\B)_+\cong \C_q[\mbox{Mat}_{s,N-s}]$. The construction of
  $\C_q[\mbox{Mat}_{N-s,s}]$ as the graded dual of $\Ubar_-$ has
  already been in detail established in \cite{a-SinShklyVa}.
\item If $\gfrak=\sofrak_{N+2}$ and $S=\pi\setminus\{\alpha_1\}$ then  
  both $(\gr\B)_-$ and $(\gr\B)_+$ are isomorphic to the
  $U_q(\sofrak_N)$-module algebra $O^N_q(\C)$ considered in
  \cite{a-FadResTak1}, the so called quantum orthogonal vector space.
\item If $\gfrak=\spfrak_{2r}$ and $S=\pi\setminus\{\alpha_{r}\}$ then  
  $(\gr\B)_-$ and $(\gr\B)_+$ are quantum coordinate algebras generated by the
  $U_q(\slfrak_r)$-module $V(2\omega_1)$ and $V(2\omega_{r-1})$,
respectively.
\item If $\gfrak=\sofrak_{2r}$, $r>3$, and $S=\pi\setminus\{\alpha_r\}$ or
  $S=\pi\setminus\{\alpha_{r-1}\}$ then 
  $(\gr\B)_-$ and $(\gr\B)_+$ are quantum coordinate algebras generated by the
  $U_q(\slfrak_r)$-module $V(\omega_2)$ and $V(\omega_{r-2})$,
respectively.
\item If $\gfrak=\efrak_6$, and $S=\pi\setminus\{\alpha_6\}$ then 
  $(\gr\B)_+$ and $(\gr\B)_-$ are quantum coordinate algebras
  generated by the spin-representation $V(\omega_5)$ of
  $U_q(\sofrak_{10})$ and its dual, respectively.
\item If $\gfrak=\efrak_7$, and $S=\pi\setminus\{\alpha_7\}$ then 
  $(\gr\B)_-$ and $(\gr\B)_+$ are quantum coordinate algebras
  generated by the $U_q(\efrak_6)$-modules $V(\omega_1)$ and $V(\omega_6)$,
  respectively.
\end{enumerate}

\section{Covariant First Order Differential Calculus}\label{diffcalc}
  For the convenience of the reader the notion of differential calculus from
\cite{a-Woro2} is recalled. 
A \textit{first order differential calculus} (FODC)
over an algebra $\B$ is a $\B$-bimodule $\Gamma$ together with a
$\C$-linear map
\begin{equation*}
  \dif:\B\rightarrow\Gamma
\end{equation*}
such that $\Gamma=\Lin_\C\{a\,\dif b\,c\,|\,a,b,c\in\B\}$ and $\dif$
satisfies the Leibniz rule
\begin{align*}
  \dif(ab)&=a\,\dif b + \dif a\,b.
\end{align*}    
Let in addition $\cA$ denote a Hopf algebra and
$\kow_\B:\B\rightarrow \cA\otimes\B$ a left $\cA$-comodule algebra
structure on $\B$.
If $\Gamma$ possesses the structure of a left $\cA$-comodule
\begin{equation*}
  \kow_\Gamma:\Gamma\rightarrow\cA\ot \Gamma
\end{equation*}
such that
\begin{equation*}
\kow_\Gamma(a\dif b\,c)=(\kow_\B a)((\id\otimes\dif)\kow_\B b)
                          (\kow_\B c)
\end{equation*}
then $\Gamma$ is called \textit{left covariant}.
For further details on first order differential calculi consult
\cite{b-KS}.

Let $U $ denote a Hopf algebra with bijective antipode and $K \subset U $
a right coideal subalgebra (For the moment $U$ and $K$ are arbitrary, in
Theorem \ref{mainTh} we will again consider the case $U=\uqg$ and
$K=\U_q(\lfrak_S)$). Consider a tensor category $\mathcal{C}$ of
finite dimensional left $U $-modules.
By this we mean as in \cite{a-MullSch} that $\cC$ is a class of finite
dimensional left $U$-modules containing the trivial
$U$-module via $\vep$ and satisfying (\ref{tenscat}). Let
$\cA :=U ^0_{\mathcal{C}}$ denote the dual Hopf algebra generated by the
matrix coefficients of all $U $-modules in $\mathcal{C}$.
Assume that $\cA$ separates the elements of $\U$ and that the antipode of
$\cA$ is bijective.
Define a left coideal subalgebra $\B \subset \cA $ by
\begin{align}\label{Bdef}
  \B :=\{b\in \cA \,|\,b_{(1)}\,b_{(2)}(k)=\vep(k)b \quad
  \text{for all $k\in K$}\}.
\end{align}
Assume that $K$ is $\mathcal{C}$-semisimple, i.e.~the restriction of any
$\U$-module in $\mathcal{C}$ to the subalgebra $K\subset \U$ is isomorphic
to the direct sum of irreducible $K$-modules. In full analogy to
\cite{a-MullSch} Theorem 2.2 (2)
this implies that $\cA $ is a faithfully flat $\B $-module.

In this situation left covariant first order differential calculi
over $\B$ can be classified via certain right ideals of
$\B^+$ \cite{a-Herm01}. More explicitly the subspace
\begin{align}
  \rid=\Big\{\sum_i \vep(a_i)b_i^+\,\Big|\,\sum_i a_i\,\dif b_i=0\Big\}
  =\{b\in\B^+\,|\,\dif b\in \B^+\Gamma\}
\end{align}
of $\B^+$, where $b^+=b-\vep(b)$ for all $b\in \B$, is a right ideal which
determines the
differential calculus uniquely.
To the FODC $\Gamma$ corresponding to this right ideal one associates the
vector space
\begin{align*}
T^\vep_\Gamma=\{f\in \B^\circ\,|\, f(x)=0 \textrm{ for all }x\in\rid\}
\end{align*}
and the so called \textit{quantum tangent space}
\begin{align*}
  T_\Gamma=(T_\Gamma^\vep)^+=\{f\in T_\Gamma^\vep\,|\, f(1)=0\}.
\end{align*}  
The dimension of a first order differential calculus is defined by
\begin{align*}
  \dim \Gamma=\dim_\C \Gamma/\B^+\Gamma=\dim_\C\B^+/\rid.
\end{align*}  

\begin{proposition}\label{corresp}{\em \cite[Cor.~5]{a-HK-QHS}}
Under the above assumptions there is a canonical one-to-one correspondence
between $n$-dimensional left covariant FODC over $\B $ and $(n+1)$-dimensional
subspaces $T^\vep\subset \B ^\circ $ such that 
\begin{align}
\vep \in T^\vep,\quad \kow T^\vep \subset  T^\vep \ot \B^\circ,\quad K T^\vep
 \subset T^\vep.
\end{align}
\end{proposition}
A covariant FODC $\Gamma\neq \{0\}$ over $\B$ is called
\textit{irreducible} if it does not possess any nontrivial quotient
(by a left covariant $\B$-bimodule). Note that this property is equivalent
to the property that $T^\vep_\Gamma$ does not possess any 
left $K$-invariant right $\B^\circ$-subcomodule $\tilde{T}$ such that
$\C\cdot\vep\varsubsetneq\tilde{T} \varsubsetneq T^\vep_\Gamma$.

We now return to the situation $U=\uqg$ and $K=\U_q(\lfrak_S)$ as in Section
\ref{U/UK+}. As an application of Theorem \ref{locfin} it is possible
to determine all finite dimensional irreducible covariant FODC over
$\B=\cqgl$ for irreducible flag manifolds $G/P_S$.

For any coalgebra $C$ let $P(C)=\{x\in C\,|\,\kopr x=1\ot x +x\ot 1\}$ denote
the vector space of primitive elements of $C$. Recall from Lemmas
\ref{coradpm} and \ref{corad} that
$P(\Ubar_+)=\Lin_\C\{E_\beta\,|\,\beta\in \gproots\}$,
$P(\Ubar_-)=\Lin_\C\{F_\beta\,|\,\beta\in \gproots\}$, and
$P(\Ubar)=P(\Ubar_+)\oplus P(\Ubar_-)$.

\begin{theorem}\label{mainTh}
  There exist exactly two nonisomorphic finite dimensional irreducible
  covariant first order differential calculi $\Gamma_+$, $\Gamma_-$
  over the quantized irreducible flag manifold $\B=\cqgl$.
  The corresponding quantum tangent spaces are
  $T_+=P(\Ubar_+)\subset \Ubar$ and $T_-=P(\Ubar_-)$.
\end{theorem}
  
\begin{proof}
  By Proposition \ref{corresp} one has to determine all finite dimensional
  left $K$-invariant right $\B^\circ$-subcomodules $T\subset F(\B^\circ,K)$
  containing $\vep$ but no other $K$-invariant $\B^\circ$-subcomodule.
  Theorem \ref{locfin} implies that $T\subset \Ubar$.
  Since $\Ubar$ is connected any $x\in C_n\Ubar$ satisfies $\kow x-1\ot
  x-x\ot 1\in \sum_{i=1}^{n-1}C_i\Ubar\ot C_{n-i}\Ubar$
  (cf.~\cite[Lem.~5.3.2(2)]{b-Montg93}). Therefore any
  nontrivial right $\Ubar$-subcomodule
  $T\subset \Ubar$ contains $\vep$ and a primitive element.
  Now the claim follows from Proposition \ref{corad} and the fact that
  $P(\Ubar_+)$ and $P(\Ubar_-)$ are irreducible $K$-modules.
\end{proof}  
  
\providecommand{\bysame}{\leavevmode\hbox to3em{\hrulefill}\thinspace}

\end{document}